# NONPARAMETRIC ESTIMATION IN A NONLINEAR COINTEGRATION TYPE MODEL

By Hans Arnfinn Karlsen, Terje Myklebust and Dag Tjøstheim

*University of Bergen, University College of Sogn og Fjordane and University of Bergen*

We derive an asymptotic theory of nonparametric estimation for a time series regression model $Z_t = f(X_t) + W_t$, where $\{X_t\}$ and $\{Z_t\}$ are observed nonstationary processes and $\{W_t\}$ is an unobserved stationary process. In econometrics, this can be interpreted as a nonlinear cointegration type relationship, but we believe that our results are of wider interest. The class of nonstationary processes allowed for $\{X_t\}$ is a subclass of the class of null recurrent Markov chains. This subclass contains random walk, unit root processes and nonlinear processes. We derive the asymptotics of a nonparametric estimate of $f(x)$ under the assumption that $\{W_t\}$ is a Markov chain satisfying some mixing conditions. The finite-sample properties of $\widehat{f}(x)$ are studied by means of simulation experiments.

**1. Introduction.** Two time series $\{X_t\}$ and $\{Z_t\}$ are said to be *linearly cointegrated* if they are both nonstationary and of unit root type and if there exists a linear combination $aX_t + bZ_t = W_t$ such that $\{W_t\}$ is stationary. This means that the series $\{X_t, Z_t\}$ move together when considered over a long period of time. The concept of cointegration was introduced by Granger [10] and further developed by Engle and Granger [6]. Since its introduction, there have been numerous papers in econometrics exploring its various aspects. Some of the main results are given in Johansen [19].

The long term relationships between two economic time series may not necessarily be linear, however, and the processes $\{X_t\}$ and $\{Z_t\}$ may not be linearly generated unit root processes. This has led to a search for nonlinear cointegration type relationships such as $Z_t = f(X_t) + W_t$, for some nonlinear function $f$ and some possibly nonlinearly generated input process $\{X_t\}$. Indeed, functional relationships of this type have been fitted to economic data (see, e.g., [8, 12]), but to our knowledge, the properties of the resulting nonparametric estimates have not been established (see [27] for a consistency property in a simplified situation, though). A brief discussion of the relationship between our work and recent contributions to the theory of nonlinear cointegration occurs in Section 6.

There are at least two difficulties (cf. [11] and others): which class of processes should be chosen as a basic class of nonstationary processes and how should an estimation theory for an estimate of $f$ be constructed? The main goal of this paper is to try to answer these questions, that is, we wish to establish a nonparametric estimation theory of the kernel estimator

$$(1.1) \qquad \widehat{f}(x) = \frac{\sum_{t=0}^n Z_t K_{x,h}(X_t)}{\sum_{t=0}^n K_{x,h}(X_t)}$$

for the function $f$ in the nonlinear regression model

$$(1.2) \qquad Z_t = f(X_t) + W_t,$$

where $K$ is a kernel function whose definition and properties are given in Section 2.1, $h$ is the bandwidth, $\{W_t\}$ is an unobserved stationary process and $\{X_t\}$ and $\{Z_t\}$ are observed processes which are nonstationary in a sense to be made precise later. At first, $\{X_t\}$ and $\{W_t\}$ will be assumed to be independent processes, which is quite a natural assumption in a nonlinear regression context. However, in a cointegration framework, this independence assumption is rather restrictive and is generally not fulfilled for linear cointegration models. In Section 4, dependence is the main subject. It turns out that dependence between $\{X_t\}$ and $\{W_t\}$ for fixed $t$ may disappear asymptotically. The reason for this phenomenon is related to restrictions on the type of dependence which is possible between a stationary and a nonstationary process. A stationary process cannot follow a nonstationary process too closely as this will violate the stationarity.

Although the connection between (1.2) and the nonlinear cointegration problem is obvious, we would like to point out that the estimation of the function $f$ in the general context we are considering should also be of interest in other areas of application. In a traditional time series regression problem, some sort of mixing condition is often assumed for $\{X_t\}$ in order to obtain a central limit theorem for $\widehat{f}(x)$. However, mixing assumptions on $\{X_t\}$ are ruled out in the general situation we consider. A minimal condition for undertaking asymptotic analysis on $\widehat{f}(x)$ is that as the number of observations on $\{X_t\}$ increases, there must be infinitely many observations in any neighborhood of $x$. This means that $\{X_t\}$ must return to a neighborhood of $x$ infinitely often, which, in turn, implies that the framework of a recurrent Markov chain is especially convenient. Since $\{X_t\}$ may be nonstationary, null recurrent processes have to be included. It should be noted that the class of null recurrent processes contains unit root processes (cf. [23]). Unlike the parametric situation, where a unit root speeds up the convergence of (global) estimates due to the large spread of the observations, in the nonparametric case, which is concerned with local estimates, the nonstationarity slows down the convergence because the time until the process returns to the local neighborhood around $x$ increases, the expected time being infinite in the null recurrent case.

In [21, 22] (hereafter, the Karlsen and Tjøstheim paper [22] is referred to as KT), an asymptotic theory was developed for nonparametric estimation for a nonstationary univariate nonlinear model in the framework of so-called $\beta$-null recurrent processes. The latter constitute a subclass of the null recurrent processes which contains the random walk. For an alternative theoretical approach in the random walk case, we refer to [26]. For a relationship between the two approaches, see [2].

We will rely on central parts of the theory of KT in our derivations in this paper. But, a host of new problems emerges in the regression case, as will be made clear in the following.

**2. Notation and some basic conditions.** We will follow the notation of KT since our proofs and results will be closely based on that paper. Thus, we denote by $\{X_t, t \geq 0\}$ a $\phi$-irreducible Markov chain on a general state space $(E, \mathcal{E})$ with transition probability $P$. This means that there exists a nontrivial measure $\phi$ on $\mathcal{E}$ such that each $\phi$-positive set $A$ is communicating with the whole state space, that is, $\sum_n P^n(x, A) > 0$ for all $x \in E$ whenever $\phi(A) > 0$, $A \in \mathcal{E}$. In this paper, we take $E \subseteq \mathbb{R}$ and we denote the class of nonnegative measurable functions with $\phi$-positive support by $\mathcal{E}^+$. For a set $A \in \mathcal{E}$, we write $A \in \mathcal{E}^+$ if the indicator function $1_A \in \mathcal{E}^+$. The process $\{X_t, t \geq 0\}$ will be assumed to be Harris recurrent. This implies that given a neighborhood $\mathcal{N}_x$ of $x$ with $\phi(\mathcal{N}_x) > 0$, $\{X_t\}$ will return to $\mathcal{N}_x$ with probability one, this being what makes asymptotics for a nonparametric estimation possible. The chain is positive recurrent if there exists an invariant probability measure such that $\{X_t, t \geq 0\}$ is strictly stationary and is null recurrent otherwise. In this paper, we are primarily interested in the null recurrent situation, in which case there exists a (unique up to a constant, nonprobability) invariant measure, which will be denoted by $\pi$.

If $\eta$ is a nonnegative measurable function and $\lambda$ is a measure, then the kernel $\eta \otimes \lambda$ is defined by

$$\eta \otimes \lambda(x, A) = \eta(x)\lambda(A), \qquad (x, A) \in (E, \mathcal{E}).$$

If $H$ is a general kernel, the function $H\eta$, the measure $\lambda H$ and the number $\lambda H\eta$ are defined, respectively, by

$$H\eta(x) = \int H(x, dy)\eta(y), \qquad \lambda H(A) = \int \lambda(dx)H(x, A),$$

$$\lambda H\eta = \int \lambda H(dy)\eta(y).$$

The convolution of two kernels, $H_1$ and $H_2$, gives another kernel, defined by

$$H_1 H_2(x, A) = \int H_1(x, dy)H_2(y, A).$$

Due to associative laws, the number $\lambda H_1 H_2 \eta$ is uniquely defined. If $A \in \mathcal{E}$ and $1_A$ is the corresponding indicator variable, then $H1_A(x) = H(x, A)$. The kernel $I_g$ is defined by $I_g(x, A) = g(x)1_A(x)$ and the special case $g = 1_C$ is denoted $I_C$.

We define $\eta \in \mathcal{E}^+$ to be *small* if there exist a measure $\lambda$, a positive constant $b$ and an integer $m \geq 1$ such that

$$(2.1) \qquad\qquad P^m \geq b\eta \otimes \lambda.$$

A set $A$ is said to be *small* if $1_A$ is small. Under quite broad conditions (cf. [9]), a compact set will be small. In this case, it follows from (2.1) that a $\phi$-positive subset of a compact set will be small. If $\lambda$ satisfies (2.1) for some $\eta$, $b$ and $m$, then $\lambda$ is a small measure.

A fundamental fact for $\phi$-irreducible Markov chains is the existence of a minorization inequality ([24], Theorem 2.1 and Proposition 2.6, pages 16–19): there exist a small function $s$, a probability measure $\nu$ and an integer $m_0 \geq 1$ such that

$$P^{m_0} \geq s \otimes \nu.$$

Some technical difficulties arise if $m_0 > 1$ because this necessitates the $m_0$-step chain; it is not a severe restriction to assume that $m_0 = 1$. Therefore, unless otherwise stated, in the sequel we will assume that the minorization inequality

$$(2.2) \qquad\qquad P \geq s \otimes \nu$$

holds, where $s$ and $\nu$ are small and $\nu(E) = 1$. In particular, this implies that $0 \leq s(x) \leq 1$, $x \in E$. If (2.2) holds, then the pair $(s, \nu)$ is called an *atom* (for $P$). A wide class of nonlinear AR(1) processes satisfying (2.2) is given in KT.

From (2.2), we obtain the identity

$$P(x, A) = (1 - s(x)) \left\{ \left( \frac{P(x, A) - s(x)\nu(A)}{1 - s(x)} \right) 1(s(x) < 1) \right.$$

$$(2.3) \qquad \qquad \left. + 1_A(x) 1(s(x) = 1) \right\} + s(x)\nu(A)$$

$$\stackrel{\text{def}}{=} (1 - s(x)) Q(x, A) + s(x)\nu(A),$$

so that the transition probability $P$ can be thought of as a mixture of the transition probability $Q$ and the small measure $\nu$. Since $\nu$ is independent of $x$, this means that the chain regenerates each time $\nu$ is chosen. This occurs with probability $s(x)$. The reasoning can be formalized by introducing the split chain $\{(X_t, Y_t)\}$, where the auxiliary chain $\{Y_t\}$ can only take values 0 and 1. Given that $X_t = x$ and $Y_{t-1} = y_{t-1}$, $Y_t$ takes the value 1 with probability $s(x)$ so that $\alpha = E \times \{1\}$ is a proper atom (cf. [24], page 51) for the split chain. We denote by

$$S_\alpha = \min\{t \geq 1 : Y_t = 1\}$$

the corresponding recurrence time. We will also make use of the consecutive sequence of recurrence times starting at time $t = 0$,

$$(2.8) \quad \tau_k = \min\{t > \tau_{k-1} : Y_t = 1\}, \qquad \tau_{-1} \stackrel{\text{def}}{=} -1 \qquad \text{for } k \geq 0, \tau = \tau_\alpha = \tau_0,$$

and the number of regenerations in the time interval $[0, n]$, that is,

$$T(n) = \max_k \{k : \tau_k \leq n\} \vee 0.$$

An invariant measure $\pi_s$ can be defined in terms of the atom $(s, \nu)$ of (2.2). In fact (KT, Section 3.2),

$$(2.9) \qquad \pi_s \stackrel{\text{def}}{=} \nu G_{s,\nu}, \qquad G_{s,\nu} \stackrel{\text{def}}{=} \sum_{\ell=0}^{\infty} (P - s \otimes \nu)^\ell.$$

If the measure $\pi_s$ is absolutely continuous with respect to Lebesgue measure, we denote by $p_s$ the corresponding density so that $p_s(x)\,dx = \pi_s(dx)$. Similarly, for $C \in \mathcal{E}^+$, we define the density $p_C(x) = p_s(x)/\pi_s 1_C$. For a $\pi_s$-integrable function $g$ on $R$, we use the notation $\pi_s g$ for

$$\pi_s g = \pi_s(g) = \int g(x)\pi_s(dx).$$

Corresponding to $T(n)$, for a set $C \in \mathcal{E}^+$, the number of times $\{X_t\}$ visits $C$ up to time $n$ is denoted by

$$T_C(n) = \sum_{t=0}^{n} 1_C(X_t).$$

From KT (Remark 3.5) we have that $T_C(n)/T(n) \xrightarrow{\text{a.s.}} \pi_s 1_C$.

The kernel $G_{s,\nu}$ of (2.9) plays an important role in Section 3 and it easily follows from the above that for a $\pi_s$-integrable $g$ defined on $E$, with $\mathbb{E}_x$ being the expectation conditional on $X(0) = x$,

$$(2.10) \qquad \mathbb{E}_x \sum_{t=0}^{\tau} g(X_t) = G_{s,\nu} g(x).$$

The minorization condition and the accompanying split chain permit the decomposition of the chain into separate and identical parts defined by the regeneration points. We have, for a function $g$,

$$(2.11) \qquad S_n(g) \stackrel{\text{def}}{=} \sum_{t=0}^{n} g(X_t) = U_0 + \sum_{k=1}^{T(n)} U_k + U_{(n)},$$

where

$$U_k = \begin{cases} \displaystyle\sum_{t=\tau_{k-1}+1}^{\tau_k} g(X_t), & \text{when } k \geq 0, \\ \displaystyle\sum_{t=\tau_{T(n)}+1}^{n} g(X_t), & \text{when } k = (n). \end{cases}$$

The sequence $\{(U_k, (\tau_k - \tau_{k-1})), \, k \geq 1\}$ consists of independent identically distributed (i.i.d.) random variables. This partition of the chain is of basic importance for the subsequent asymptotic analysis. In the following, we will sometimes use the symbol $U = U(g)$ to denote a random variable representing the common marginal distribution of $\{U_k, \, k \geq 1\}$.

We must introduce a restriction on the way the process regenerates: the chain $\{X_t\}$ is $\beta$-null recurrent if there exist a small nonnegative function $f$, an initial measure $\lambda$, a constant $\beta \in (0, 1)$ and a slowly varying function $L_f$ such that

$$(2.12) \qquad \mathbb{E}_\lambda \sum_{t=0}^{n} f(X_t) \sim \frac{1}{\Gamma(1+\beta)} n^\beta L_f(n)$$

as $n \to \infty$. This condition is equivalent to (cf. KT, Theorem 3.1) a restriction on the tail distribution of the recurrence time $S_\alpha$, in that

$$(2.13) \qquad \mathbb{P}_\alpha(S_\alpha > n) = \frac{1}{\Gamma(1-\beta) n^\beta L_s(n)} (1 + \mathcal{O}(1)),$$

where $L_s$ is a slowly varying function depending on $s$ and where $\mathbb{P}_\alpha$ means that the initial distribution is equal to $\delta_\alpha(x, y)$, that is, $Y_0 = 1$, $X_0 = x$ arbitrary. In the sequel, (2.13) will be referred to as the *tail condition*.

A random walk process is $\beta$-null recurrent with $\beta = 1/2$.

2.1. *Basic conditions.* We denote by $h = h_n$ the bandwidth used in the nonparametric estimation. It is assumed to satisfy $h_n \to 0$ and, with no loss of generality, we also assume that $h_n \leq 1$. Let $K : \mathbb{R} \to \mathbb{R}$ be a kernel function and for a fixed $x$, let $K_{x,h}(y) = h^{-1}K((y-x)/h)$, $\mathcal{N}_x(h) = \{y : K_{x,h}(y) \neq 0\}$ and $\mathcal{N}_x = \mathcal{N}_x(1)$. In our context, a locally bounded function will be taken to mean a function bounded in a neighborhood of $x$ and a locally continuous function is a function continuous at the point $x$. Without loss of generality, we may assume that this neighborhood equals $\mathcal{N}_x$ and that local continuity implies local boundedness. This follows since $\mathcal{N}_x(h) = x \oplus h\mathcal{N}_0$.

We will consider the problem of evaluating the properties of the kernel estimator (1.1) of the function $f$ of (1.2) under the assumption that $\{W_t\}$ is Markov. In Section 3, $\{X_t\}$ and $\{W_t\}$ are assumed to be independent. The independence assumption is removed in Section 4, and the compound process $\{(X_t, W_t)\}$ is assumed to be Markov.

The following set of conditions is always assumed:

$B_0$  (i) the kernel $K$ is nonnegative, $\int K(u)\, du < \infty$ and
$\qquad\quad$ $\|K\|_2^2 = \int K^2(u)\, du < \infty$;
$\qquad$(ii) the $\{X_t\}$ process is a Harris recurrent Markov chain;
$\qquad$(iii) the transfer function $f$ is continuous at the point $x$.

We will also make heavy use of the following conditions $B_1$–$B_4$ of KT. For ease of reference, these conditions are restated here:

$B_1$      (i) $\int K(u)\, du = 1$;
$\qquad$  (ii) $\int uK(u)\, du = 0$;
$B_2$      (i) the support $\mathcal{N}_0$ of the kernel is contained in a compact set;
$\qquad$  (ii) the kernel is bounded and $\mathcal{N}_x$ is a small set;
$B_3$  the invariant measure $\pi_s$ has a locally continuous density $p_s$ which is
$\qquad$ locally strictly positive, that is, $p_s(x) > 0$;
$B_4$  for all $\{A_h\} \in \mathcal{E}$ such that $A_h \downarrow \varnothing$, $\lim_{h \downarrow 0} \overline{\lim}_{y \to x} P(y, A_h) = 0$.

In all of the proofs, we use $c_1, c_2, \ldots$ as a sequence of generic constants and if $\{a_n\}$ and $\{b_n\}$ are two real-valued strictly positive sequences, then we write $a_n \ll b_n$ if $a_n = \mathcal{O}(b_n)$. The associated $\sigma$-algebra, $\mathcal{F}_t^X$, for a stochastic process $\{X_t,\ t \geq 0\}$ is defined in the usual way: $\mathcal{F}_t^X = \sigma\{X_j, j \leq t\}$ and $\mathcal{F}^X = \bigvee_t \mathcal{F}_t^X$.

**3. Nonparametric estimation of $f$.** At the outset, we assume that $\{W_t\}$ is a $\phi$-irreducible ergodic Markov chain which satisfies (2.2). Additional assumptions will be introduced as needed. Actually, we also allow a slight generalization of (1.2), in that we include an instantaneous transformation of $W_t$, resulting in

$$(3.1) \qquad\qquad Z_t = f(X_t) + g_W(W_t).$$

This is an extension of (1.2) since even if $\{W_t\}$ is Markov, $\{W_t'\} = \{g_W(W_t)\}$ does not have to be Markov. We assume that $\mathbb{E}g_W(W_0) = 0$. Because we do not generally restrict $g_W$ to be a small function (consider, e.g., $g_W(w) \equiv w$), Lemma 5.1 and Lemma 5.2 of KT cannot be used, which complicates matters considerably.

Throughout this section, we make the assumption that $\{X_t\}$ and $\{W_t\}$ are independent and, using this assumption, we are able to obtain results which are of interest in the general context of nonparametric estimation of nonstationary processes. In Section 4, we allow for dependence, but put restrictions on $g_W$, and some parts of the results obtained in this section are extended. Moreover, our findings in Section 4 highlight the fact that the actual dependence occurring in cointegration models disappears asymptotically. In this way, results in this section are also relevant to cointegration models. Furthermore, they may serve as a starting point for deriving asymptotic results for the dependent case without the restrictions on $g_W$ which are imposed in Section 4. We believe that letting $W_t' = g_W(X_t, \ldots, X_{t-p}, W_t)$ for some fixed $p$, where $\{W_t\}$ is a Markov process, independent of $\{X_t\}$ and such that $\{W_t'\}$ is stationary, may be a possible way to proceed.

We start by expressing $\hat{f}(x) - f(x)$ in the $S_n$-notation of (2.11), and this is done by rewriting the numerator of $\hat{f}(x)$ of (1.1) as

$$Z_t = g_W(W_t) + (f(X_t) - f(x)) + f(x),$$

$$Z_t K_{x,h}(X_t) = g_h(X_t, W_t) + \psi_x(X_t)K_{x,h}(X_t) + f(x)K_{x,h}(X_t),$$

where $g_h(z, u) = g_W(u) \cdot K_{x,h}(z)$ and $\psi_x(y) = f(y) - f(x)$. By the definition of $\hat{f}(x)$, this gives

$$\hat{f}(x) - f(x) = S_n^{-1}(K_{x,h})\{S_n(g_h) + S_n(\psi_x \cdot K_{x,h})\}.$$

The last term on the right-hand side represents the bias. It is a stochastic quantity and we want to replace it by a deterministic bias term. Let

$$a_h \overset{\text{def}}{=} \frac{\pi_s I_{K_{x,h}} \psi_x}{\pi_s I_{K_{x,h}} 1}, \qquad b_h \overset{\text{def}}{=} I_{K_{x,h}}(\psi_x - a_h).$$

Then

$$\begin{aligned}
\hat{f}(x) - f(x) - a_h &= S_n^{-1}(K_{x,h})\{S_n(g_h) + S_n(\psi_x \cdot K_{x,h}) - a_h S_n(K_{x,h})\} \\
&= S_n^{-1}(K_{x,h})\{S_n(g_h) + S_n(b_h)\}.
\end{aligned}$$

It follows that

$$(3.2) \qquad h^{1/2} S_n^{1/2}(K_{x,h})\{\hat{f}(x) - f(x) - a_h\} = \Delta_{n,h}^1 + \Delta_{n,h}^2,$$

where

$$\Delta^1_{n,h} = S_n^{-1/2}(K_{x,h})h^{1/2}S_n(g_h),$$

$$\Delta^2_{n,h} = \{\widehat{p}_C(x)\}^{-1/2}T_C^{-1/2}(n)h^{1/2}S_n(b_h),$$

$$\widehat{p}_C(x) = T_C^{-1}(n)S_n(K_{x,h})$$

and where $C$ is a purely auxiliary small set. Replacing $P\xi$ with $f$ in the proof of Theorem 5.4 in KT, then using $B_1$–$B_3$ and condition (2.13), we have $\Delta^2_{n,h_n} = \mathcal{O}_P(1)$ and by KT (the second part of Theorem 5.3 and also the proof of Theorem 5.4), $\widehat{p}_C(x) = p_C(x) + \mathcal{O}_P(1)$ since $\pi_s b_h = 0$.

By (3.2), the above arguments show that a central limit theorem for $\widehat{f}(x)$ follows from a central limit theorem for $\Delta^1_{n,h_n}$. We continue the proof of the asymptotic properties of $\widehat{f}$ by formulating a general nonparametric CLT.

### 3.1. A nonparametric CLT for null recurrent processes.

Assume that $\{X_t\}$ is a general Markov chain [e.g., it could be identified with the compound chain $\{(X_t, W_t)\}$ or with just one of the components] which satisfies the minorization condition (2.2) and the tail condition (2.13). Let (assuming first- and second-order moments exist)

$$U_0 = U_0(g_h) = \sum_{t=0}^{\tau_0} g_h(X_t), \qquad \mu(g_h) = \mathbb{E}U(g_h),$$

$$\sigma^2(g_h) = \mathbb{E}U^2(g_h) - \mu^2(g_h),$$

where $g_h$ is a real-valued function defined on $E$ for all $h > 0$ and $\tau = \tau_0$ is defined as in (2.8). Note that with the function $g_h$ used in this paper, the random variables $U_1(g_h), U_2(g_h), \ldots$ in the decomposition (2.11) are independent so that in the notation of equation (4.4) of KT, $\overline{\sigma}^2(g_h) = \sigma^2(g_h)$. Consider the following conditions, where $-\infty < \mu, \mu' < \infty$, $0 < \sigma, \sigma' < \infty$, $v \in [0,1]$, $m \geq 2$, $\epsilon > 0$, $0 < d_m, d'_m < \infty$, $\beta$ is defined in (2.13), $\lambda$ is an initial measure and $h \downarrow 0$.

C1: $\mu(g_h) = \mu + \mathcal{O}(1)$, $\mu(|g_h|) = \mu' + \mathcal{O}(1)$.
C2: $h\sigma^2(g_h) = \sigma^2 + \mathcal{O}(1)$.
C3: $h\sigma^2(|g_h|) = \sigma'^2 + \mathcal{O}(1)$.
C4: $\mathbb{E}|U(g_h) - \mu(g_h)|^{2m} \leq d_m h^{-2m+v}$.
C5: $\mathbb{E}|U(|g_h|) - \mu(|g_h|)|^{2m} \leq d'_m h^{-2m+v}$.
C6: $h_n^{-1} \ll n^{\beta\delta_m - \epsilon}$, $\delta_m = \frac{m-1}{m-v}$.
C7: $\exists g_0 : h|g_h| \leq g_0$ and $\mathbb{P}_\lambda(U_0(g_0) < \infty) = 1$.

The following theorem is essentially a translation of a CLT result in KT. It will be used to prove the main CLT results of the present paper.

THEOREM 3.1. *Let $C$ be a small set. Assume that the tail condition* (2.13) *and* C1–C6 *hold with $\mu = 0$ for an $m \geq 2$ and a $v \in [0, 1]$. Then for any initial measure $\lambda$ for $X_0$ such that* C7 *holds,*

$$h_n^{1/2} T_C^{-1/2}(n)\{S_n(g_{h_n}) - T_C(n)\pi_s^{-1}(C)\mu(g_{h_n})\} \xrightarrow[n]{d} \mathcal{N}(0, \sigma^2 \pi_s^{-1}(C)).$$

PROOF. The proof is essentially based on KT (Theorem 4.2). Since $g_h$ is a function of one variable, the conditions in that theorem simplify. Clearly, conditions $A_0$–$A_2$ of Theorem 4.2 of KT follow directly from C1–C3. In conditions C4 and C5 the quantity $v$ is allowed to vary everywhere in $[0, 1]$, whereas in conditions $A_3$ and $A_4$ of KT, $v$ can only take the values 0 and 1. However, this extension is allowed by a trivial modification of the first part of the proof of Theorem 4.1 of KT. Condition $A_5$ of Theorem 4.2 of KT follows straightforwardly from C7 by reasoning as in the proof of Theorems 5.1 and 5.3 of KT.  □

Before we can employ Theorem 3.1, we need to analyze the regeneration structure of $\{(X_t, W_t)\}$ more carefully. This is done in a series of lemmas in Sections 3.2–3.8. We believe that these results are of independent interest and that they are potentially useful in other situations. Our main result is stated in Section 3.9.

3.2. *Decomposition of $S_n(g)$.* We assume that the compound chain $\{(X_t, W_t)\}$ satisfies (2.2) so that it can be extended by the split chain method, with $\{(X_t, W_t, Y_t)\}$ being a split chain. Note that if $\{X_t\}$ and $\{W_t\}$ separately satisfy the minorization inequality (2.2), it is not obvious that the compound chain $\{(X_t, W_t)\}$ will. However, if $\{X_t\}$ and $\{W_t\}$ are independent, then it is trivial to verify (2.2), as is shown at the beginning of Section 3.3. Let

$$\tau_k = \inf\{t > \tau_{k-1} : Y_t = 1\}, \qquad k \geq 0, \tau_{-1} = -1.$$

Then the sequence $\{\tau_k\}$ represents the regeneration times for the compound process. The basic decomposition, (2.11), with $g = g_h$ defined at the beginning of this section, gives

$$(3.3) \quad S_n(g) = U_0(g) + \sum_{k=1}^{T(n)} U_k(g) + U_{(n)}(g), \qquad T(n) = \sup\{k : \tau_k \leq n\} \vee 0,$$

where

$$U_k(g) = \begin{cases} \displaystyle\sum_{t=\tau_{k-1}+1}^{\tau_k} g(X_t, W_t), & \text{for } k \geq 0, \\ \displaystyle\sum_{t=\tau_{T(n)}+1}^{n} g(X_t, W_t), & \text{for } k = (n). \end{cases}$$

According to the general theory, the variables $\{(U_k(g), (\tau_k - \tau_{k-1})), \, k \geq 1\}$ are i.i.d. We denote by $U = U(g)$ a random variable having the common marginal distribution of the $U_k$'s and write $\mu(g) = \mathbb{E}U = \mathbb{E}_\nu U_0(g)$, $\sigma^2(g) = \mathrm{Var}(U) = \mathrm{Var}_\nu(U_0(g)) = \mathbb{E}_\nu U_0^2(g) - \mu^2(g)$, where $\nu$ refers to the compound chain $\{(X_t, W_t)\}$.

Our first problem is to find conditions which ensure that $\mu(|g|)$ and $\sigma^2(g)$ are finite. Again, by reference to the general theory [cf. Appendix A, (A.11) and (A.12)], we have that, with $s$ referring to the compound chain,

$$(3.4) \qquad \mu(g) = \pi_s g, \qquad \sigma^2(g) = \pi_s g^2 + 2\pi_s I_g H \, G_{s,\nu} g - \pi_s^2 g,$$

where $H = P - s \otimes \nu$ and $G_{s,\nu}$ is defined as in (2.9).

The conditions ensuring $\sigma^2(g) < \infty$ are not evident from (3.4) if we want to avoid the relatively strong restriction that $g_W$ is a small function. If $g_W(w) \equiv w$, then requiring $g_W$ to be small is roughly equivalent to $\phi$-mixing, which is not satisfied for, say, an autoregressive process. The problem is linked to the term $G_{s,\nu}$. In fact, we also need to demonstrate the existence of higher moments and to verify conditions connected to the bandwidth as seen in C1–C7.

### 3.3. $\beta$-null recurrence for the compound process.

Let $P$ denote the transition probability for the Markov process $\{(X_t, W_t)\}$. We label quantities associated with $\{X_t\}$ by 1 and with $\{W_t\}$ by 2. The transition probability $P$ satisfies (2.2) when $P_1$ and $P_2$ do since

$$(3.5) \qquad P = P_1 \otimes P_2 \geq (s_1 \otimes s_2) \otimes (\nu_1 \otimes \nu_2) = s \otimes \nu.$$

Condition (3.5) will be assumed to hold in the following.

LEMMA 3.1. *Assume that $\{X_t\}$ and $\{W_t\}$ are independent, that the tail condition (2.13) holds for $\{X_t\}$ and that $\{W_t\}$ is ergodic. Then the compound process $\{(X_t, W_t)\}$ is $\beta$-null recurrent, that is, the tail condition holds for the compound process.*

PROOF. Let $C_1$ and $C_2$ be small sets and let $\nu = \nu_1 \otimes \nu_2$. Then

$$(3.6) \qquad \begin{aligned} \mathbb{E}_\nu \left\{ \sum_{t=0}^{n} 1_{C_1}(X_t) 1_{C_2}(W_t) \right\} &= \sum_{t=0}^{n} (\nu_1 P_1^t 1_{C_1})(\nu_2 P_2^t 1_{C_2}) \\ &= (\pi_2 1_{C_2}) \sum_{t=0}^{n} \nu_1 P_1^t 1_{C_1} + \sum_{t=0}^{n} (\nu_1 P_1^t 1_{C_1}) b_t, \end{aligned}$$

where $b_t = \nu_2 P_2^t 1_{C_2} - \pi_2 1_{C_2}$ and where $\pi_2$ is the stationary measure for $\{W_t\}$. Since $\{W_t\}$ is ergodic, $b_t = \mathcal{O}(1)$. Since $\{X_t\}$ is $\beta$-null [cf. KT, Lemma 3.1 and formulas (3.12) and (3.13)], we have that

$$\sum_{t=0}^{n} \nu_1 P_1^t 1_{C_1} = (\pi_{s_1} 1_{C_1}) \psi_1(n)(1 + a_n), \qquad \psi_1(n) = n^\beta L_{s_1}(n), \,\, a_n = \mathcal{O}(1).$$

By (2.12), the conclusion of the lemma follows if we can show that the second term of (3.6) is $\mathcal{O}(\psi_1(n))$. Let $\psi_M = \sup_{t \le M} \psi_1(t)$, $A = \sup_t |a_t|$, $c_1 = \pi_{s_1} 1_{C_1}$, $B = \sup_t |b_t|$ and $B^{(M)} = \sup_{t > M} |b_t|$. Then for all $M > 0$,

$$(3.7) \qquad \frac{\sum_{t=0}^n (\nu_1 P_1^t 1_{C_1}) |b_t|}{\psi_1(n)} \le c_1 \left\{ B \frac{\psi_M (1+A)}{\psi_1(n)} + B^{(M)} (1 + |a_n|) \right\}.$$

Letting $n$ tend to infinity and then letting $M$ tend to infinity, we find that the left-hand side of (3.7) is $\mathcal{O}(1)$ with respect to $n$. $\square$

3.4. *Refinement of the decomposition structure.* We extend both chains with the split chain method and write $\{(X_t, Y_t^1)\}$ and $\{(W_t, Y_t^2)\}$. Due to independence, $\{(X_t, W_t, Y_t)\}$ is the split chain for the compound process $\{(X_t, W_t)\}$, where $Y_t = Y_t^1 Y_t^2$ (cf. [24], (4.17), page 62). Thus,

$$(3.8) \qquad \tau_k = \inf\{t > \tau_{k-1} : Y_t^1 = Y_t^2 = 1\}, \qquad k \ge 0, \tau_{-1} = -1.$$

We shall now look more closely at the decomposition structure and try, to some extent, to reduce it to the marginal decomposition of the $\{X_t\}$-process, that is, the regenerations defined by $\{\tau_k^1\}$,

$$(3.9) \qquad \tau_k^1 = \inf\{t > \tau_{k-1}^1 : Y_t^1 = 1\}, \qquad k \ge 0, \tau_{-1}^1 = -1,$$

which defines the $X$-partition. Let

$$(3.10) \qquad V_j = V_j(g) = \sum_{t=\tau_{j-1}^1 + 1}^{\tau_j^1} g(X_t, W_t), \qquad s \ge 0.$$

Although the $V_j$'s are neither unconditionally nor conditionally independent, they will be useful. By (3.8), we see that the regeneration times for the compound chain are also regeneration times for $\{X_t\}$. Hence, following the regeneration times (3.9) for the $X$-process, we recover all of the simultaneous regeneration times given by (3.8). The gaps between successive simultaneous regeneration times define a subdivision of each $U_k$ into $V_j$'s and this refines the decomposition given by (3.3). Let

$$\mathcal{T}_k = \inf\{j > \mathcal{T}_{k-1} : Y_{\tau_j^1}^2 = 1\} \qquad \text{for } k \ge 0, \mathcal{T}_{-1} = -1, \mathcal{T} \overset{\text{def}}{=} \mathcal{T}_0.$$

Then $\tau_k = \tau_{\mathcal{T}_k}^1$, which gives

$$U_k = \sum_{t=\tau_{\mathcal{T}_{k-1}}^1 + 1}^{\tau_{\mathcal{T}_k}^1} g(X_t, W_t) = \sum_{j=\mathcal{T}_{k-1}+1}^{\mathcal{T}_k} \sum_{t=\tau_{j-1}^1 + 1}^{\tau_j^1} g(X_t, W_t) = \sum_{j=\mathcal{T}_{k-1}+1}^{\mathcal{T}_k} V_j$$

and in particular for $k = 0$,

$$(3.11) \qquad U_0 = \sum_{j=0}^{\mathcal{T}} V_j.$$

The number of subblocks inside a large block is distributed as the recurrence time for the ergodic process $\{W_{\tau_k^1}\}$. Comparing this distribution with $\tau$ and $\tau^1$, it is evident that a block which is quite large is partitioned into relatively few sub-blocks. The advantage of this construction is that the subblocks are defined by the regeneration times for the $X$-process and the $X$-part of $g_h$ is marginally a small function.

3.5. *The embedded process.* The following lemma proves that embedding the $\{X_t\}$-regeneration times into $\{W_t\}$ and extending to a split chain are essentially commutative operations.

LEMMA 3.2. *The process $\{W_{\tau_k^1}, k \geq 0\}$ is a Markov process with transition probability $\underset{\sim}{P} = P_2 \Phi_{\nu_1}$, where $\Phi_{\nu_1} = \sum_{\ell=0}^{\infty} \{\nu_1 (P_1 - s_1 \otimes \nu_1)^{\ell} s_1\} P_2^{\ell}$. Moreover,*

$$\underset{\sim}{P} \geq \underset{\sim}{s} \otimes \underset{\sim}{\nu}, \tag{3.12}$$

*with $(\underset{\sim}{s}, \underset{\sim}{\nu}) = (s_2, \nu_2 \Phi_{\nu_1})$. Let $\lambda = \lambda_1 \otimes \lambda_2$ be the initial measure for $\{(X_t, W_t)\}$. Let $\{\widehat{W}_k\} = \{(W_k, Y_k)\}$ be the split chain generated by $\underset{\sim}{P}$ and $(\underset{\sim}{s}, \underset{\sim}{\nu})$ and let $\{\widehat{W}_{\tau_k^1}\} = \{(W_{\tau_k^1}, Y_{\tau_k^1}^2)\}$. Then*

$$\{\widehat{W}_{\tau_k^1}\} \overset{\mathrm{d}}{=} \{\widehat{\underset{\sim}{W}}_k\} \tag{3.13}$$

*when the initial measure for $\underset{\sim}{W}_0$ is $\underset{\sim}{\lambda} = \widetilde{\lambda} \overset{\mathrm{def}}{=} \lambda_2 \Phi_{\lambda_1}$. In particular, let $\underset{\sim}{\mathcal{T}}$ denote the first regeneration time for $\{\widehat{\underset{\sim}{W}}_k\}$. Then the occupation time formula is given by*

$$\mathbb{E}_{\lambda} \sum_{k=0}^{\mathcal{T}} 1_A(W_{\tau_k^1}) = \mathbb{E}_{\widetilde{\lambda}} \sum_{k=0}^{\underset{\sim}{\mathcal{T}}} 1_A(\underset{\sim}{W}_k) = \begin{cases} \widetilde{\lambda} \underset{\sim}{G}_{\underset{\sim}{s}, \underset{\sim}{\nu}} 1_A, & \text{in general,} \\ \pi_2 \underset{\sim}{G}_{\underset{\sim}{s}, \underset{\sim}{\nu}} 1_A, & \text{if } \lambda_2 = \pi_2, \\ \pi_{s_2} 1_A, & \text{if } \lambda = \nu, \end{cases} \tag{3.14}$$

*where $\underset{\sim}{G}_{\underset{\sim}{s}, \underset{\sim}{\nu}} = \sum_{\ell=0}^{\infty} (\underset{\sim}{P} - \underset{\sim}{s} \otimes \underset{\sim}{\nu})^{\ell}$.*

The proof is given in Appendix B.

Intuitively, changing the time parameter from $\{k\}$ to $\{\tau_k^1\}$ in the ergodic process $\{W_t\}$ should decrease the amount of dependence, and this is the content of our next result. More specifically, we obtain that the rate of convergence of the transition probability toward the stationary measure is at least as good for the $\{\underset{\sim}{W}_k\}$-process as for the $\{W_t\}$-process.

LEMMA 3.3. *Suppose that $\{W_t\}$ is geometric ergodic. Then this is also true for $\{\underset{\sim}{W}_k\}$. If $\{W_t\}$ is strongly mixing with mixing rate defined by $\alpha =$* $\{\alpha_j\}$, then $\{\underline{W}_k\}$ is strongly mixing with mixing rate $\underline{\alpha}$, which is equal to or faster than $\alpha$. In particular, for an integer $p \geq 0$,

$$(3.15) \qquad \sum_{\ell=1}^{\infty} \ell^p \alpha_\ell < \infty \quad \Longrightarrow \quad \mathbb{E}_{\pi_2} \underline{\mathcal{T}}^{p+1} < \infty.$$

The proof is given in Appendix B.

3.6. *Moment bounds.* Our nonparametric CLT requires bounds for the moments of $U(g)$ given by C4 and C5. We first need to find upper bounds for moments of $U(g)$ corresponding to (3.11) and related quantities. We assume that

$$(3.16) \qquad g(x, w) = (g_X \otimes g_W)(x, w) = g_X(x)g_W(w).$$

Our method is to use a representation of $U(g)$ as a partial sum of $V$'s, these variables being defined by the regeneration of $\{X_t\}$.

In the following, $H_j = P_j - s_j \otimes \nu_j$ for $j = 1, 2$ and as before, $H = P - s \otimes \nu$. Also, recall that $I_g$ is defined by $I_g(x, A) = g(x)1_A(x)$.

THEOREM 3.2. *Let $m \geq 1$ and $V_j$ be defined by* (3.10) *and* (3.16). *Then*

$$(3.17) \qquad \mathbb{E}_\nu \sum_{j=0}^{\mathcal{T}} |V_j|^m \leq \pi_{s_2} |g_W|^m \mathbb{E} U^m(|g_X|).$$

*For all $p > 0$ and $\delta \in (0, \infty)$,*

$$(3.18) \quad \mathbb{E}|U(g)|^p \leq \mathbb{E}_{\nu'}^{1/(1+\delta)} \left\{ \sum_{j=0}^{\mathcal{T}} |V_j|^{p(1+\delta)} \right\} \mathbb{E}_{\nu'}^{\delta/(1+\delta)} |\mathcal{T} + 1|^{p(1+\delta^{-1})}.$$

The proof will be based on two lemmas. We use the notation $\delta_j = \tau_j^1 - \tau_{j-1}^1$ for $j \geq 0$ and $\mathcal{H}_j = \mathcal{F}_{\tau_j^1}^X \vee \mathcal{F}_{\tau_j^1}^{Y^1} \vee \mathcal{F}^W \vee \mathcal{F}^{Y^2}$. Then $V_j$ is measurable $\mathcal{H}_j$, and $\{\mathcal{T}_0 \geq j\} = \{\mathcal{T} \geq j\} \in \mathcal{H}_{j-1}$. By (3.11), $U_0 = \sum_{j=0}^{\infty} V_j 1(\mathcal{T} \geq j)$ and for $m \geq 1$,

$$\mathbb{E}_\lambda \{V_j^m 1(\mathcal{T} \geq j)\} = \mathbb{E}_\lambda \{1(\mathcal{T} \geq j)\mathbb{E}_\lambda[V_j^m \mid \mathcal{H}_{j-1}]\}.$$

The following technical result, which is the first step in the proof of Theorem 3.2, uses the independence of $\{X_t\}$ and $\{W_t\}$, together with the regeneration property of $\{X_t\}$.

LEMMA 3.4 [Decoupling]. *Let $\lambda = \lambda_1 \otimes \lambda_2$. Let $j \geq 0$ be fixed and let $\{X_t'\}$ be an independent copy of $\{X_t\}$ so that $\{X_t'\}$ is independent of both* $\{X_t\}$ and $\{W_t\}$. Let $\xi_W$ be a real-valued function defined on $R \times \{0, 1\}$ and for fixed $j$, let

$$a_\ell^j = \xi_W(W_{\tau_{j-1}^1 + \ell + 1}, Y_{\tau_{j-1}^2}^2), \qquad \ell \geq 0, Y_{\tau_{-1}^2}^2 = y$$

and let $V_{\xi,j}$ be an extension of (3.10), given by

$$(3.19) \qquad V_{\xi,j} = \sum_{t=\tau_{j-1}^1 + 1}^{\tau_j^1} g_X(X_t)\xi_W(W_t, Y_{\tau_{j-1}^1}^2), \qquad j \geq 0.$$

Then for $m \geq 1$,

$$\mathbb{E}_\lambda\{V_{\xi,j}^m \mid \mathcal{H}_{j-1}\} = \begin{cases} \mathbb{E}_{\lambda_1} U_0^m(\underline{a}, g_X), & \text{for } j = 0, \\ \mathbb{E}_{\nu_1} U_0^m(\underline{a}, g_X), & \text{for } j \geq 1, \end{cases}$$

where $U_0(\underline{a}, g_X) = \sum_{\ell=0}^{\tau_0^1} g_X(X_\ell')a_\ell$ and $\underline{a} = \{a_\ell\} = \{a_\ell^1\}$.

PROOF. Let $j \geq 1$. By (3.19),

$$V_{\xi,j} = \sum_{t=\tau_{j-1}^1+1}^{\tau_j^1} g_X(X_t)\xi_W(W_t, Y_{\tau_{j-1}^2}^2) = \sum_{\ell=1}^{\tau_j^1 - \tau_{j-1}^1} g_X(X_{\tau_{j-1}^1 + \ell})\xi_W(W_{\tau_{j-1}^1 + \ell}, Y_{\tau_{j-1}^2}^2),$$

so that with $S_\alpha^1 = \min\{t \geq 1 : Y_k^1 = 1\}$,

$$\mathbb{E}_\lambda\{V_{\xi,j}^m \mid \mathcal{H}_{j-1}\} = \mathbb{E}_\lambda\left\{\left[\sum_{\ell=1}^{\delta_j} g_X(X_{\tau_{j-1}^1 + \ell})\xi_W(W_{\tau_{j-1}^1 + \ell}, Y_{\tau_{j-1}^2}^2)\right]^m \Big| \mathcal{H}_{j-1}\right\}$$

$$= \mathbb{E}_\alpha\left\{\sum_{\ell=1}^{S_\alpha^1} g_X(X_\ell')a_{\ell-1}\right\}^m$$

$$= \mathbb{E}_{\nu_1}\left\{\sum_{\ell=0}^{\tau_0^1} g_X(X_\ell')a_\ell\right\}^m$$

$$= \mathbb{E}_{\nu_1} U_0^m(\underline{a}, g_X),$$

where we have used the fact that

$$\mathcal{L}_{\lambda_1}\{X_{\tau_{j-1}^1 + \ell}, 1 \leq \ell \leq \delta_j\} = \mathcal{L}_\alpha\{X_\ell', 1 \leq \ell \leq S_\alpha^1\} = \mathcal{L}_{\nu_1}\{X_\ell', 0 \leq \ell \leq \tau_0^1\},$$

with $\mathcal{L}_\lambda$ denoting the simultaneous distribution with initial measure $\lambda$.

If $j = 0$, then

$$\mathbb{E}_\lambda\{V_{\xi,j}^m \mid \mathcal{H}_{j-1}\} = \mathbb{E}_\lambda\left\{\left[\sum_{t=0}^{\tau_0^1} g_X(X_t)\xi_W(W_t, y)\right]^m \Big| \mathcal{F}^W \vee \mathcal{F}^{Y^2}\right\}$$

$$= \mathbb{E}_{\lambda_1} U_0^m(\underline{a}, g_X).$$

$\square$

Using the previous lemma, the factorization of $g$ given by (3.16) and a general moment formula given in Corollary A.1 in Appendix A, we obtain a useful exact formula. The notation is in accordance with Theorem A.1 in Appendix A. We use the index set $\Delta_r^m = \{\alpha \in \mathcal{N}_+^r : \sum_{j=1}^r \alpha_j = m\}$, where $\mathcal{N}_+^r$ is the $r$-Cartesian product of all strictly positive integers, the multinomial coefficient $\binom{m}{\alpha} = \frac{m!}{\alpha_1! \cdots \alpha_r!}$, $\mathcal{N}_{0,+}^r = \mathcal{N} \times \mathcal{N}_+^{r-1}$ and $j_{(2)} = (j_2, \ldots, j_r) \in \mathcal{N}_+^{r-1}$.

LEMMA 3.5.   Let $V_k$ be defined by (3.10). Let $m \geq 1$. Then

$$(3.20) \quad \mathbb{E}_\nu \sum_{k=0}^{\mathcal{T}} V_k^m = \sum_{r=1}^m \sum_{\alpha \in \Delta_r^m} \binom{m}{\alpha} \sum_{j_{(2)} \in \mathcal{N}_+^{r-1}} \{\pi_{s_1} \widehat{f}_{j_{(2)},\alpha}^X\}\{\pi_{s_2} \widehat{f}_{j_{(2)},\alpha}^W\},$$

where

$$\widehat{f}_{j_{(2)},\alpha}^X = I_{g_X^{\alpha_1}} H_1^{j_2} I_{g_X^{\alpha_2}} \cdots H_1^{j_r} I_{g_X^{\alpha_r}} 1,$$

$$\widehat{f}_{j_{(2)},\alpha}^W = I_{g_W^{\alpha_1}} P_2^{j_2} I_{g_W^{\alpha_2}} \cdots P_2^{j_r} I_{g_W^{\alpha_r}} 1.$$

More generally, we have for $\lambda = \lambda_1 \otimes \lambda_2$, with $f_{j,\alpha}^X = H_1^{j_1} \widehat{f}_{j_{(2)},\alpha}^X$ and $f_{j,\alpha}^W = P_2^{j_1} \widehat{f}_{j_{(2)},\alpha}^W$,

$$(3.21) \quad \begin{aligned} \mathbb{E}_\lambda \sum_{k=0}^{\mathcal{T}} V_k^m &= \sum_{r=1}^m \sum_{\alpha \in \Delta_r^m} \binom{m}{\alpha} \sum_{j \in \mathcal{N}_{0,+}^r} (\lambda - \nu)(f_{j,\alpha}^X \otimes f_{j,\alpha}^W) \\ &\quad + \{\nu_1 f_{j,\alpha}^X\}\{\widetilde{\lambda} \underset{\sim s, \nu}{G} P_2 f_{j,\alpha}^W\}. \end{aligned}$$

REMARK 3.1.   If $\lambda = \nu$, then

$$\widetilde{\lambda} \underset{\sim s, \nu}{G} P_2 f_{j,\alpha}^W = \underset{\sim}{\nu} \underset{\sim s, \nu}{G} P_2 P_2^{j_1} \widehat{f}_{j_{(2)},\alpha}^W = \pi_{s_2} P_2^{j_1+1} \widehat{f}_{j_{(2)},\alpha}^W = \pi_{s_2} \widehat{f}_{j_{(2)},\alpha}^W,$$

$$\sum_{j_1=0}^\infty \nu_1 f_{j,\alpha}^X = \nu_1 \sum_{j_1=0}^\infty H_1^{j_1} \widehat{f}_{j_{(2)},\alpha}^X = \pi_{s_1} \widehat{f}_{j_{(2)},\alpha}^X.$$

Thus, (3.21) reduces to (3.20) when $\lambda = \nu$.

PROOF OF LEMMA 3.5.   We rewrite the first term on the left-hand side of (3.20) using the fact that $1(\mathcal{T} \geq k) = 1(\mathcal{T} \geq k-1)1(Y_{\tau_{k-1}^1}^2 = 0)$, so that

$$\mathbb{E}_\nu \sum_{k=0}^{\mathcal{T}} V_k^m = \mathbb{E}_\nu(V_0^m) + \sum_{k=1}^\infty \mathbb{E}_\nu \{1(\mathcal{T} \geq k-1)\mathbb{E}_\nu(V_k^m 1(Y_{\tau_{k-1}^1}^2 = 0) \mid \mathcal{H}_{k-1})\}.$$

Let $V_{\xi,k}$ be defined by (3.19), where $\xi_W(w,y) = g_W(w)(1-y)$. By Lemma 3.4 and its proof, it is seen that the conditional mean given $\mathcal{H}_{k-1}$ only involves the regeneration of $\{X_t\}$ and we can therefore use Appendix A (and more specifically Corollary A.1) to obtain for $k \geq 1$,

$$\mathbb{E}_\nu[V_k^m 1(Y_{\tau_{k-1}^2}^2 = 0) \mid \mathcal{H}_{k-1}]$$

$$= \mathbb{E}_\nu[V_{\xi,k}^m \mid \mathcal{H}_{k-1}]$$

$$= \sum_{r=1}^m \sum_{\alpha \in \Delta_r^m} \binom{m}{\alpha} \sum_{j \in \mathcal{N}_{0,+}^r} \{\nu_1 f_{j,\alpha}^X\} \left[\prod_{i=1}^r g_W^{\alpha_i}(W_{\tau_{k-1}+t_i+1})\right] 1(Y_{\tau_{k-1}^1}^2 = 0),$$

where $t_i = j_1 + \cdots + j_i$. Let $\mathcal{G}_k = \mathcal{F}_{\tau_k^1}^X \vee \mathcal{F}_{\tau_k^1}^{Y^1} \vee \mathcal{F}_{\tau_k^1}^W \vee \mathcal{F}_{\tau_{k-1}^1}^{Y^2}$. Then by conditioning with respect to $\mathcal{G}_{k-1}$, we find that

$$\mathbb{E}_\nu \left\{ 1(\mathcal{T} \geq k-1) \left[\prod_{i=1}^r g_W^{\alpha_i}(W_{\tau_{k-1}^1+t_i+1})\right] 1(Y_{\tau_{k-1}^1}^2 = 0)\right\}$$

$$= \mathbb{E}_\nu\{1(\mathcal{T} \geq k-1) H_2 f_{j,\alpha}^W(W_{\tau_{k-1}^1})\}.$$

Hence,

$$\mathbb{E}_\nu \sum_{k=1}^{\mathcal{T}} V_k^m = \sum_{k=1}^\infty \mathbb{E}_\nu\{1(\mathcal{T} \geq k-1)\mathbb{E}_\nu[V_k^m 1(Y_{\tau_{k-1}^2}^2 = 0) \mid \mathcal{H}_{k-1}]\}$$

$$(3.22) \qquad = \sum_{r=1}^m \sum_{\alpha \in \Delta_r^m} \binom{m}{\alpha} \sum_{j \in \mathcal{N}_{0,+}^r} \{\nu_1 f_{j,\alpha}^X\} \mathbb{E}_\nu \left\{\sum_{k=1}^\infty 1(\mathcal{T} \geq k-1) H_2 f_{j,\alpha}^W(W_{\tau_{k-1}^1})\right\}$$

$$= \sum_{r=1}^m \sum_{\alpha \in \Delta_r^m} \binom{m}{\alpha} \sum_{j \in \mathcal{N}_{0,+}^r} \{\nu_1 f_{j,\alpha}^X\} \{\pi_{s_2} H_2 f_{j,\alpha}^W\}.$$

Similarly, we find that

$$(3.23) \qquad \mathbb{E}_\nu V_0^m = \sum_{r=1}^m \sum_{\alpha \in \Delta_r^m} \binom{m}{\alpha} \sum_{j \in \mathcal{N}_{0,+}^r} \{\nu_1 f_{j,\alpha}^X\} \{\nu_2 f_{j,\alpha}^W\}$$

and by combining (3.22) and (3.23) and using $\pi_{s_2} H_2 = \pi_{s_2} - \nu_2$, we get

$$\mathbb{E}_\nu \sum_{k=0}^{\mathcal{T}} V_k^m = \sum_{r=1}^m \sum_{\alpha \in \Delta_r^m} \binom{m}{\alpha} \sum_{j \in \mathcal{N}_{0,+}^r} \{\nu_1 f_{j,\alpha}^X\} \{\pi_{s_2} f_{j,\alpha}^W\}.$$

Now, $\pi_{s_2} f_{j,\alpha}^W = \pi_{s_2} \widehat{f}_{j(2),\alpha}^W$ and

$$\sum_{j_1=0}^\infty \nu_1 f_{j,\alpha}^X = \sum_{j_1=0}^\infty \nu_1 H_1^{j_1} I_{g_X^{\alpha_1}} H_1^{j_2} I_{g_X^{\alpha_2}} \cdots H_1^{j_r} I_{g_X^{\alpha_r}} 1 = \pi_{s_1} \widehat{f}_{j(2),\alpha}^X$$

and thus (3.20) is proved.

The proof of (3.21) is similar. Instead of (3.22), we obtain

$$
\begin{aligned}
\mathbb{E}_\lambda \sum_{k=1}^{\mathcal{T}} V_k^m &= \sum_{r=1}^m \sum_{\alpha \in \Delta_r^m} \binom{m}{\alpha} \sum_{j \in \mathcal{N}_{0,+}^r} \{\nu_1 f_{j,\alpha}^X\} \mathbb{E}_\lambda \sum_{k=0}^{\mathcal{T}} H_2 f_{j,\alpha}^W(W_{\tau_k^1}) \\
&= \sum_{r=1}^m \sum_{\alpha \in \Delta_r^m} \binom{m}{\alpha} \sum_{j \in \mathcal{N}_{0,+}^r} \{\nu_1 f_{j,\alpha}^X\} \{\widetilde{\lambda} G_{\underset{\sim}{g}, \underline{\nu}} H_2 f_{j,\alpha}^W\}
\end{aligned}
\tag{3.24}
$$

and (3.23) is changed to

$$
\mathbb{E}_\lambda V_0^m = \sum_{r=1}^m \sum_{\alpha \in \Delta_r^m} \binom{m}{\alpha} \sum_{j \in \mathcal{N}_{0,+}^r} \{\lambda_1 f_{j,\alpha}^X\} \{\lambda_2 f_{j,\alpha}^W\}.
\tag{3.25}
$$

Combining (3.24) and (3.25) and using $\widetilde{\lambda} G_{\underset{\sim}{g}, \underline{\nu}} H_2 = \widetilde{\lambda} G_{\underset{\sim}{g}, \underline{\nu}} P_2 - \nu_2$, we obtain (3.21). $\square$

REMARK 3.2. If $m = 1$, then by (3.20), $\mu(g) = \mathbb{E}_\nu \sum_{j=0}^{\mathcal{T}} V_j = \{\pi_{s_1} g_X\} \{\pi_{s_2} g_W\}$, which, using (3.11), is consistent with (3.4).

REMARK 3.3. If $m = 2$, then

$$
\mathbb{E}_\nu \sum_{j=0}^{\mathcal{T}} V_j^2 = \{\pi_{s_1} g_X^2\} \{\pi_{s_2} g_W^2\} + 2 \sum_{\ell=1}^\infty \{\pi_{s_1} I_{g_X} H_1^\ell I_{g_X} 1\} \{\pi_{s_2} I_{g_W} P^\ell I_{g_W} 1\}.
$$

REMARK 3.4. By (3.11) and (3.21), we find that for general $\lambda = \lambda_1 \otimes \lambda_2$, $g = g_X \otimes g_W$, we have

$$
\mathbb{E}_\lambda U_0(g) = \sum_{j=0}^\infty (\lambda - \nu)(P_1^j g_X \otimes P_2^j g_W) + \sum_{j=0}^\infty \{\nu_1 H_1^j g_X\} \{\widetilde{\lambda} G_{\underset{\sim}{g}, \underline{\nu}} P_2^{j+1} g_W\}.
$$

If $g_X$ is small, $\lambda_2 = \pi_2$ and $\sup_j \pi_2 G_{\underset{\sim}{g}, \underline{\nu}} P_2^{j+1} |g_W| < \infty$, then $\mathbb{E}_\lambda U_0(g)$ is finite. More generally, taking $p = 1 + \delta$, $f = P_2^{j+1} |g_W|$ and $\lambda = \lambda_2 = \pi_2$ in Lemma A.2, we have that

$$
\pi_2 G_{\underset{\sim}{g}, \underline{\nu}}^{(1+\delta)/(1+\eta\delta)} P_2^{j+1} |g_W| \le c_2 \mathbb{E}_{\pi_2}^{1/(1+\eta\delta)} \underset{\sim}{\mathcal{T}}^{1+2\delta} \{\pi_2^{\eta\delta/(1+\eta\delta)} |g_W|^{(1+\delta)/(\eta\delta)}\},
$$

with $\eta \in (0,1)$ and $\delta > 0$ arbitrary. The result can now be combined with Lemma 3.3, which ensures the existence of moments of $\underset{\sim}{\mathcal{T}}$ under appropriate mixing conditions.

PROOF OF THEOREM 3.2. We first prove (3.17). Assume that $g_X \geq 0$. By Cauchy–Schwarz, recalling that $\sum_{j=1}^r \alpha_j = m$, we have

$$
\begin{aligned}
(3.26) \qquad |\pi_2 \widehat{f}_{j(2),\alpha}^W| &= |\pi_2 I_{g_W^{\alpha_1}} P_2^{j_1} I_{g_W^{\alpha_2}} \cdots P_2^{j_{r-1}} I_{g_W^{\alpha_r}} 1| \\
&\leq \prod_{r=1}^m \pi_2^{[\alpha_r/m]} |g_W|^{[m/\alpha_r]\alpha_r} \leq \pi_2 |g_W|^m.
\end{aligned}
$$

Inserting (3.26) into (3.20), we obtain

$$
\mathbb{E}_\nu \sum_{k=0}^{\mathcal{T}} V_k^m \leq \pi_{s_2} |g_W|^m \sum_{r=1}^m \sum_{\alpha \in \Delta_r^m} \binom{m}{\alpha} \sum_{j(2) \in \mathcal{N}_+^{r-1}} \pi_{s_1} I_{g_X^{\alpha_1}} H_1^{j_2} I_{g_X^{\alpha_2}} \cdots H_1^{j_r} I_{g_X^{\alpha_r}} 1
$$

$$
= \pi_{s_2} |g_W|^m \mathbb{E} U^m(g_X),
$$

from which (3.17) follows trivially.

To prove (3.18), let $r = 1 + \delta$ and $q = 1 + \delta^{-1}$. Then $\mathbb{E}_\nu |U_0(g)|^p = \mathbb{E}_\nu |\sum_{j=0}^{\mathcal{T}} V_j|^p$ and

$$
\begin{aligned}
\mathbb{E}_\nu \left| \sum_{j=0}^{\mathcal{T}} V_j \right|^p &\leq \mathbb{E}_\nu \max_{0 \leq j \leq \mathcal{T}} |V_j|^p |\mathcal{T} + 1|^p \leq \mathbb{E}_\nu^{1/r} \max_{0 \leq j \leq \mathcal{T}} |V_j|^{pr} \mathbb{E}_\nu^{1/q} |\mathcal{T} + 1|^{pq} \\
&\leq \mathbb{E}_\nu^{1/r} \sum_{j=0}^{\mathcal{T}} |V_j|^{pr} \mathbb{E}_\nu^{1/q} |\mathcal{T} + 1|^{pq}. \qquad \square
\end{aligned}
$$

3.7. *Moment bounds of $U(g_h)$ expressed in terms of bandwidth.* The following results describe how higher-order moments of $U$ behave as functions of the bandwidth. This is what is needed to apply C4 and C5 in Theorem 3.1.

THEOREM 3.3. *Let $g_X = g_{X,h} = K_{x,h}$ and assume that conditions B$_2$, B$_3$ and (3.16) hold. Then for all integers $k, m \geq 1$,*

$$
(3.27) \qquad \mathbb{E}|U(g_h)|^{2m} \leq d_{m,k} h^{-2m+1/(k+1)},
$$

*where*

$$
d_{m,k} \stackrel{\text{def}}{=} \{\pi_{s_2}^{1/(k+1)} g_W^{2m(k+1)}\} \mathbb{E}_\nu^{k/(k+1)} |\mathcal{T} + 1|^{2m(k+1/k)} \{d_{2m}'^{1/(k+1)}\}
$$

*and the sequence of constants $\{d_m'\}$ is only dependent on $\mathcal{N}_x$ and $\sup_u K(u)$.*

PROOF. By (3.18), with $p = 2m$ and $\delta = k$, we have

$$
(3.28) \quad \mathbb{E} U^{2m}(g_h) \leq \mathbb{E}_\nu^{1/(k+1)} \left| \sum_{j=0}^{\mathcal{T}} |V_j|^{2m(k+1)} \right| \mathbb{E}_\nu^{k/(k+1)} |\mathcal{T} + 1|^{2m((k+1)/k)}.
$$

From (3.17), we have

$$(3.29) \qquad \mathbb{E}_\nu \sum_{j=0}^{\mathcal{T}} V_j^{2m(k+1)} \le \pi_{s_2} |g_W|^{2m(k+1)} \mathbb{E} U^{2m(k+1)}(|g_{X,h}|),$$

and by KT [Lemma 5.2 with $\xi_0 \equiv 1$ and $\mathcal{G}_2$ replaced by $\mathcal{G}_1 = \{f : (E, \mathcal{E}) \mapsto (\mathbb{R}, \mathcal{B}(\mathbb{R}))\}$, where $\mathcal{B}(\mathbb{R})$ is the class of all Borel sets on $\mathbb{R}$],

$$(3.30) \qquad \mathbb{E} U^{2m(k+1)}(|g_{X,h}|) \le d'_{2m} h^{-2m(k+1)+1}.$$

In the proof of that lemma, it is also shown that the sequence of constants $\{d'_m\}$ is only dependent on $\mathcal{N}_x$ and $\sup_u K(u)$.

Inserting (3.29) and (3.30) into (3.28), we get (3.27). □

3.8. *Asymptotic variance.* Exact information about the first order properties of the asymptotic variance is important (cf. C2 and C3). Such information is contained in the next result, which is the analogue of Lemma 5.1 of KT. Our method of proof uses a truncation technique based on the notion of a generalized autocovariance function. We believe the latter concept to be of some independent interest.

THEOREM 3.4. *Assume that the process $\{W_t\}$ is an irreducible, ergodic, strongly $\alpha$-mixing process which satisfies (2.2) and has mixing rate satisfying $\sum_\ell \ell^{[2/k] \vee 1} \alpha_\ell < \infty$, $\pi_2 g_W = 0$ and $\pi_2 |g_W|^{2(k+1)} < \infty$ for some integer $k \ge 1$.*

*Assume, in addition, that $g_{X,h} = K_{x,h}$ and that conditions B2–B4 hold. Then if $\mu(g_W) = 0$, we have, as $h \downarrow 0$,*

(i) $h\sigma^2(g_{X,h} \otimes g_W) = p_{s_1}(x) \|K\|_2^2 \pi_{s_2} g_W^2 + \mathcal{O}(1)$,
(ii) $h\sigma^2(|g_{X,h} \otimes g_W|) = h\sigma^2(g_{X,h} \otimes g_W) + \mathcal{O}(1)$.

In the proof of Theorem 3.4, we need some results of a more general nature concerning generalized autocovariances. These are formulated in Lemmas 3.6 and 3.7 below, for a general $\phi$-irreducible, aperiodic, Harris recurrent Markov chain $\{X_t\}$ with transition function $P$ satisfying (2.2) and with the $S_n(g)$-decomposition, as in (2.11). The next step is Lemma 3.8, where we apply Lemmas 3.6 and 3.7 to our Markov chain $\{(X_t, W_t)\}$ [with a slight conflict of notation, taking $X_t = (X_t, W_t)$].

We begin by extending the notion of a cross-covariance function, as defined for ergodic processes.

DEFINITION 3.1. *Let $g, f \in L^1(\pi_s) \cap L^2(\pi_s)$. The generalized covariance and cross-covariance function is defined by*

$$(3.31) \qquad \gamma_{g,f}(\ell) = \begin{cases} \pi_s I_{g_0} f_0 + \mu_g \mu_f (1 - \pi_s s^2), & \text{when } \ell = 0, \\ \varphi_g P^{\ell-1} f_0, & \text{when } \ell \ge 1, \\ \gamma_{f,g}(-\ell), & \text{when } \ell < 0, \end{cases}$$

where $\varphi_g \stackrel{\text{def}}{=} \pi_s I_g P - \mu_g \nu$ and $\gamma_g \stackrel{\text{def}}{=} \gamma_{g,g}$. Mean centering a function $f$ with $s\mu_f$ produces a function denoted by $f_0$, that is, $f_0 \stackrel{\text{def}}{=} f - s\mu_f$.

Note that when $\mu_g = \mu_f = 0$, the generalized cross-covariance function is equal to the ordinary one, apart from the constant $c_\pi = \pi s$, that is, $\gamma_{g,f} = c_\pi^{-1}\overline{\gamma}_{g,f}$, where $\overline{\gamma}_{g,f}$ denotes the stationary covariance function.

LEMMA 3.6. *Assume that $|g|$ is small. Then $\sigma^2(g) = \sum_{\ell=-\infty}^{\infty} \gamma_g(\ell)$.*

PROOF. We have, by (A.12) in Appendix A and by (3.4), that

$$\sigma^2(g) = \gamma_g(0) + 2\varphi_g G_{s,\nu} g_0.$$

Iterating $G_{s,\nu} = I + (P - s \otimes \nu)G_{s,\nu}$, we get $G_{s,\nu} = G^{(n)} + P^n G_{s,\nu} - G^{(n)} s \otimes \pi_s$ with $G^{(n)} = \sum_{\ell=0}^{n-1} P^\ell$. Pre- and post-multiplying this equation by $\varphi_g$ and $g_0$, respectively, gives $\varphi_g G_{s,\nu} g_0 = \sum_{\ell=1}^{n} \gamma_g(\ell) + \varphi_g P^n \psi$, where $\psi \stackrel{\text{def}}{=} G_{s,\nu} g_0$. By Nummelin ([24], Theorem 6.7, page 109), and since $\|\varphi_g\| \leq 2\mu_{|g|}$, $|g_0|$ is small and $\psi$ is bounded, we find that $\varphi_g P^n \psi = \mathcal{O}(1)$, from which the result follows. □

REMARK 3.5. The formula in Lemma 3.6 can be viewed as a generalization of the formula $\text{Var}(n^{-1/2}\sum_{j=0}^{n} X_j) = \sum_{\ell=-\infty}^{\infty} \text{Cov}(X_t, X_{t-\ell}) + \mathcal{O}(1)$ in the case where $\{X_t\}$ is a stationary process with an absolutely summable covariance function.

It is necessary to weaken the assumption of smallness in Lemma 3.6.

LEMMA 3.7. *Assume that*

(i) $g \in L^1(\pi_s) \cap L^2(\pi_s)$, $\pi_s I_{|g|} P G_{s,\nu} |g| < \infty$.

*If there is an approximating sequence $\{g_n\}$, in the sense that $|g_n|$ is a small function, $|g_n| \leq |g|$ and $g_n \xrightarrow[n]{} g$ a.s. $[\pi_s]$, then for each $\ell$,*

(ii) $\gamma_{g_n}(\ell) \xrightarrow[n]{} \gamma_g(\ell)$;

(iii) $\sigma^2(g) = \lim_n \sum_{\ell=-\infty}^{\infty} \gamma_{g_n}(\ell)$.

*Suppose that*

(iv) $\gamma_{g_n}(\ell) = \eta_n b_\ell + d_{n,\ell}$,

*where $\sum_{\ell=1}^{\infty} \sup_n |d_{n,\ell}| < \infty$, $\eta_n = \eta + \mathcal{O}(1)$ and $\sum_{\ell=1}^{\infty} b_\ell < \infty$. Then*

(v) $\sigma^2(g) = \sum_{\ell=-\infty}^{\infty} \gamma_g(\ell)$, *and if $\eta = 0$, the convergence is absolute.*

PROOF. Let $\{g_n\}$ be an approximating sequence which satisfies the conditions in the lemma. First, we prove that

$$(3.32) \qquad \lim_n \pi_s I_{g_n} PG_{s,\nu} g_n = \pi_s I_g PG_{s,\nu} g.$$

Let $\xi_n = PG_{s,\nu} g_n$, $\xi = PG_{s,\nu} g$ and $\xi_0 = PG_{s,\nu}|g|$ so that $|\xi_n| \leq \xi_0$. We must show that $\xi_n \xrightarrow{n} \xi = PG_{s,\nu} g$ a.s. $[\pi_s]$. Let $D$ be the set of points where $g_n$ fails to converge toward $g$. Then $\pi_s\{G1_D > 0\} = 0$ with $G = \sum_{\ell=0}^{\infty} P^\ell$ since $\pi_s$ is a maximal irreducible measure. Hence, $\pi_s\{PG_{s,\nu}1_D > 0\} = 0$. The rest of the proof of (3.32) follows directly from the dominated convergence theorem since $\pi_s I_{g_n} PG_{s,\nu} g_n = \pi_s I_{g_n} \xi_n = \pi_s(g_n \xi_n)$, where $|g_n \cdot \xi_n| \leq \{|g| \cdot \xi\}$ and $(g_n \cdot \xi) \xrightarrow{n} (g \cdot \xi)$ a.s. $[\pi_s]$.

By Lemma 3.6 and (3.32), statement (iii) holds. It is obvious that (ii) holds and if (iv) is true, then $\sum_{\ell=1}^{\infty} \gamma_{g_n}(\ell) \xrightarrow{n} \sum_{\ell=1}^{\infty} \gamma_g(\ell)$, by the dominated convergence theorem. Together with (iii), we can conclude that (v) is true. □

In the next lemma, we return to the Markov chain $\{(X_t, W_t)\}$ and let it play the role of the general Markov chain in Lemmas 3.6 and 3.7.

LEMMA 3.8. *Assume that $g(x,w) = g_1(x)g_2(w)$, $g_1$ is small and that $\{W_t\}$ satisfies the conditions stated in Theorem 3.4. Then $\sigma^2(g) = \sum_{\ell=-\infty}^{\infty} \gamma_g^{X,W}(\ell)$.*

PROOF. Our proof is based on Lemma 3.7. We must show that Lemma 3.7(i) and Lemma 3.7(iv) are satisfied. We do not assume that $\mu_{g_2} = 0$.

By (3.28), (3.29) with $m = 1$ and the smallness of $g_1$, we have

$$(3.33) \qquad \pi_s(g^2) + 2\pi_s I_g PG_{s,\nu} g \leq \mathbb{E}U^2(|g|)$$
$$\leq c_0 \pi_2^{1/(k+1)} |g_2|^{2(k+1)} \mathbb{E}_\nu^{k/(k+1)} \mathcal{T}^{(2k+2)/k}.$$

The quantity $\pi_2 |g_2|^{2(k+1)}$ is finite by assumption. We have that $\mathbb{E}_\nu \mathcal{T}^{(2k+2)/k} = \mathbb{E}_{\pi_2} \mathcal{T}^{(2k+2)/k}$, thus the right-hand side of (3.33) is finite if $\mathbb{E}_{\pi_2} \mathcal{T}^{(k+2)/k} < \infty$ (cf. [4, 5]). By Lemma 3.3, this is true if $\sum_{\ell=1}^{\infty} \ell^{[2/k]\vee 1} \underline{\alpha}_\ell < \infty$ and is thus satisfied by the mixing assumption on $\{W_t\}$. Hence, Lemma 3.7(i) holds.

We begin the next step by establishing an approximating sequence for $g = g_1 \otimes g_2$. By Nummelin ([24], Corollary 2.1, page 24), there exists an increasing sequence of small sets $C'_n$ such that $\bigcup_{n=1}^{\infty} C'_n = E_2$, where $E = E_1 \times E_2$. Define $C_n = (|g_2| < n) \cap C'_n$. Let

$$(3.34) \qquad g_n = g_1 \otimes g_2^n, \qquad g_2^n = g_2 1_{C_n}.$$

Then $|g_n|$ is small for all $n$ and $|g_n| \uparrow |g|$ a.s. $[\pi_s]$.

We express $\gamma_{g_1 \otimes g_2^n}$ in terms of $g_1$ and $g_2^n$ using Definition 3.1. At the same time, we insert $\Delta_2^\ell = (\nu_2 P_2^\ell - \pi_2)$, $c_\pi = \pi_2 s$, $\overline{\gamma}_{g_2} = c_\pi \gamma_{g_2}$ and $\overline{\mu}_{g_2} = c_\pi \mu_{g_2}$. This gives, for $\ell > 0$,

$$
\begin{aligned}
\gamma_{g_n}(\ell) = & \{\pi_{s_1} I_{g_1} P_1^\ell g_1\}\{\pi_{s_2} I_{g_2} P_2^\ell g_{20}^n\} + \mu_{g_2^n}\{\varphi_{g_1} P_1^{\ell-1} g_1\}\{\Delta_2^{\ell-1} g_{20}^n\} \\
(3.35) \qquad & + \mu_{g_2^n}\{\pi_{s_1} I_{g_1} P_1^\ell g_{10}\}\{\overline{\gamma}_{g_2^n}(\ell) - \overline{\mu}_{g_2^n}\Delta_2^{\ell-1} s_2\} \\
& + \mu_{g_2^n}^2 \gamma_{g_1}(\ell)\{\Delta_2^{\ell-1} s_2 + \pi_2 s_2\}.
\end{aligned}
$$

By the mixing property of $\{W_t\}$, we find that $\overline{\gamma}_{g_2}$ is absolutely summable (cf. [14], Corollary A.2, page 278). Moreover, since the recurrence time for $\{W_t\}$ has a finite second-order moment, $\{W_t\}$ is ergodic of degree 3 as a Markov chain (cf. [24], page 84) and that implies the finiteness of $\sum_{\ell=1}^\infty \ell \|\Delta_2^\ell\|$ (cf. [24], Theorem 6.13, page 118). By Lemma 3.6, $\gamma_{g_1}$ is summable. It is now easy to verify that each of the four terms of $\gamma_{g_n}(\ell)$ given by (3.35) satisfies Lemma 3.7(iv). Hence, by Lemma 3.7(v), the proof is finished. $\square$

PROOF OF THEOREM 3.4. By Lemma 3.8 and Definition 3.1, since $\mu_{g_W} = 0$, we have

$$
\sigma^2(g_{X,h} \otimes g_W) = \sum_\ell \{\pi_{s_1} I_{K_{x,h}} P_1^\ell K_{x,h}\} \gamma_{g_W}(\ell).
$$

For $\ell > 0$, $h\pi_{s_1} I_{K_{x,h}} P_1^\ell K_{x,h} = \mathcal{O}(1)$, by $B_4$ [cf. KT, proof of part (b) of Lemma 5.1]. Since

$$
|\pi_{s_1} I_{K_{x,h}} P_1^\ell K_{x,h}| \leq \int p_s(x+hu)K(u)P^\ell(x+hu, \mathcal{N}_x(h)) \, du \leq c_0
$$

and $\sum |\gamma_{g_W}(\ell)|$ is finite, we can apply the dominated convergence theorem, that is,

$$
\begin{aligned}
(3.36) \qquad \lim_{h\downarrow 0} h\sigma^2(g_{X,h} \otimes g_W) &= \sum_\ell \lim_{h\downarrow 0}\{h\,\pi_{s_1} I_{K_{x,h}} P_1^\ell K_{x,h}\} \gamma_{g_W}(\ell) \\
&= \lim_{h\downarrow 0}\{h\,\pi_{s_1} K_{x,h}^2\} \gamma_{g_W}(0) = p_{s_1}(x)\|K\|_2^2 \pi_{s_2} g_W^2.
\end{aligned}
$$

The proof of Theorem 3.4(ii) follows in a similar way, using Lemma 3.8 and (3.35) with $g_2^n = g_W$. For $\ell = 0$, Definition 3.1 must be used. $\square$

3.9. *Main result.*

THEOREM 3.5. *Assume that $\{X_t\}$ and $\{W_t\}$ are independent, recurrent Markov chains in (3.1) and that $\{W_t\}$ is an irreducible, ergodic, strongly $\alpha$-mixing process which satisfies (2.2) and has a mixing rate satisfying $\sum_\ell \ell^{[2/k]\vee 1}\alpha_\ell < \infty$, $\pi_2 g_W = 0$ and $\pi_2 |g_W|^{2m(k+1)} < \infty$ for some integers $k \geq 1$* and $m \geq 2$, $\pi_2$ being the invariant probability measure of $\{W_t\}$. Moreover, assume that $\mathrm{B_1}$–$\mathrm{B_4}$ hold and that (2.2) and the tail condition (2.13) hold for $\{X_t\}$.

Finally, assume that for some $\epsilon > 0$,

$$h_n^{-1} \ll n^{\beta \delta_m - \epsilon}, \qquad \delta_m = \frac{m-1}{m - 1/(k+1)}.$$

Then for all $\lambda = \lambda_1 \otimes \pi_2$, we have

$$\left\{ h_n \sum_{t=0}^{n} K_{x,h_n}(X_t) \right\}^{1/2} \left\{ \widehat{f}(x) - f(x) - \frac{\pi_s I_{K_{x,h_n}} \psi_x}{\pi_s K_{x,h_n}} \right\} \xrightarrow[n]{d} \mathcal{N}(0, \sigma_W^2 \|K\|_2^2).$$

If the density $p_s$ and the function $f$ possess continuous derivatives of second order, then the bias term $\pi_s I_{K_{x,h_n}} \psi_x / \pi_s K_{x,h_n}$ is negligible when $h_n^{-1} \gg n^{\beta/5 + \epsilon}$.

PROOF. We use Theorem 3.1 on the compound chain $\{(X_t, W_t)\}$. As noted at the beginning of Section 3, it is enough to prove that

$$\Delta_{n,h_n} = S_n^{-1/2}(K_{x,h_n}) h_n^{1/2} S_n(K_{x,h_n} \otimes g_W) \xrightarrow[n]{d} \mathcal{N}(0, \|K\|_2^2 \pi_2 g_W^2).$$

Recall that for $C = C_1 \times C_2$, $C_i \in \mathcal{E}_i$, $i = 1, 2$,

$$T_C(n) = \sum_{t=0}^{n} 1_{C_1}(X_t) 1_{C_2}(W_t)$$

represents the number of visits of $\{(X_t, W_t)\}$ to $C$ up to time $n$. We choose $C_1$ and $C_2$ so that both sets are small. Then by KT (the second part of Theorem 5.3), using $\mathrm{B_2}$–$\mathrm{B_4}$ and the tail condition (2.13), we have

$$\widehat{p}_{C_1}(x) \overset{\text{def}}{=} \frac{S_n(K_{x,h_n})}{T_{C_1 \times E_2}(n)} = p_{C_1}(x) + \mathcal{O}_P(1),$$

with $E = E_1 \times E_2$ and where $p_{C_1}(x) = p_{s_1}(x) / \pi_{s_1} 1_{C_1}$. By KT (Remark 3.5),

$$\frac{T_C(n)}{T_{C_1 \times E_2}(n)} = \frac{\pi_s 1_C}{\pi_s 1_{C_1 \times E_2}} + \mathcal{O}(1) = \pi_2(C_2) + \mathcal{O}(1) \qquad \text{a.s.}$$

We can write

$$\Delta_{n,h} = \widehat{p}_{C_1}^{-1/2}(x) \left\{ \frac{T_C(n)}{T_{C_1 \times E_2}(n)} \right\}^{1/2} \{ T_C^{-1/2}(n) h^{1/2} S_n(g_h) \}$$

$$= A_{n,h}^{1/2} \Delta_{n,h}^0,$$

say, where $g_h(z, w) = K_{x,h}(z) g_W(w)$ and $A_{n,h_n} = \{ p_{C_1}^{-1}(x) \pi_2 1_{C_2} \} + \mathcal{O}_P(1)$. Hence, it is enough to prove that

$$(3.37) \qquad \Delta_{n,h_n}^0 \xrightarrow[n]{d} \mathcal{N}\left( 0, p_{C_1}(x) \|K\|_2^2 \frac{\pi_2 g_W^2}{\pi_2 1_{C_2}} \right).$$

By $\mathbf{B_3}$ and Bochner's theorem, $\mathbf{C1}$ is satisfied. From Theorem 3.3 and Theorem 3.4, conditions $\mathbf{C2}$–$\mathbf{C5}$ are satisfied with $v = 1/(k+1)$.

It only remains to verify $\mathbf{C7}$. Let $g_0 = c_0 1_{\mathcal{N}_x} |g_W|$, where $c_0$ is an appropriate constant. Then $|h g_h| \le g_0$. We must prove that $\mathbb{P}_\lambda(U_0(g_0) < \infty) = 1$, with [cf. (3.11)] $U(g_0) = \sum_{j=0}^{\mathcal{T}} V_j(g_0)$. But, this is satisfied if $\mathbb{E}_\lambda U_0(g_0) < \infty$. By Remark 3.4, this is true if

$$\mathbb{E}_{\pi_2} |\mathcal{T}|^{1+2\delta} |\pi_2| g_W |^{(1+\delta)/(\eta\delta)} < \infty \tag{3.38}$$

for some $\delta > 0$ and $\eta \in (0,1)$. Let $k \ge 1$ be fixed and $\pi_2 |g_W|^{2(k+1)} < \infty$. By Lemma 3.3, (3.38) is satisfied if

$$\frac{1+\delta}{\eta\delta} \le 2(k+1), \qquad 1 + 2\delta \le 2 + \frac{2}{k}.$$

This is true if

$$\frac{1}{2k+1} < \delta < 1 + \frac{1}{k}, \qquad \frac{1+\delta}{\delta(2k+1)} \le \eta < 1.$$

Thus, (3.38) holds.

Hence, by Theorem 3.1, $\Delta_{n,h_n}^0 \xrightarrow{d}{n} \mathcal{N}(0, \sigma_C^2)$ and by Theorem 3.4,

$$\sigma_C^2 = \{\pi_{s_1} 1_{C_1}\}^{-1} \{\pi_{s_2} 1_{C_2}\}^{-1} p_{s_1}(x) \|K\|_2^2 \pi_{s_2} g_W^2.$$

It follows that (3.37) holds. $\square$

REMARK 3.6. If $k = 1$, then we require that the residual $\{W_t\}$-process have a finite eighth-order moment, together with a mixing rate which satisfies $\sum_\ell \ell^2 \alpha_\ell < \infty$. If, on the other hand, all moments of the residual process are finite, then it is enough for there to exist a $\delta > 0$ such that

$$\sum_\ell \alpha_\ell^{1-\delta} < \infty.$$

## 4. Some extensions to the dependent case.

In linear cointegration theory, the stationary process $\{W_t\}$ resulting from a linear cointegration relationship $W_t = Z_t - a X_t$, say, will generally be dependent on $\{X_t\}$. From this point of view, it is of interest to extend the theory of Section 3. We will do this by assuming that $\{(X_t, W_t)\}$ is a Markov chain in (3.1) and specifying a dependence relation between them for which the asymptotic theory holds. In this situation, we will prove that the compound process $\{(X_t, W_t)\}$ is $\beta$-null recurrent, as was done in the previous section. But, unlike Section 3, we essentially assume that the function $(u, w) \mapsto K_{x,h}(u) g_W(u, w)$ is small. In this way, it is guaranteed that the necessary moment requirements are satisfied. In addition, we need existence and smoothness of an invariant measure for the compound chain, together with additional conditions which control the bias.

4.1. *Conditional expectation.* The restriction on the type of dependence allowed between $\{X_t\}$ and $\{W_t\}$ will be formulated in terms of the conditional expectation of $W_t$ with respect to $X_t$. Let $\{(X_t, W_t)\}$ be Harris null recurrent with state space $(E, \mathcal{E}) = (E_1 \times E_2, \mathcal{E}_1 \otimes \mathcal{E}_2)$, invariant measure $\pi_s$ and maximal irreducibility measure $\phi$. Assume that

$$\pi_s 1_{C_1 \times E_2} < \infty \qquad \text{for some } C_1 \in \mathcal{E}_1$$

and let

$$(4.1) \qquad Q \overset{\text{def}}{=} \frac{\pi_s I_{C_1 \times E_2}}{\pi_s 1_{C_1 \times E_2}},$$

so that $Q$ is a probability measure on $(E', \mathcal{E}') = (E_1' \times E_2, \mathcal{E}_1' \otimes \mathcal{E}_2)$ with $E_1' = C_1$ and $\mathcal{E}_1' = \mathcal{E}_1 \cap C_1$. Here, $I_{C_1 \times E_2}$ is defined as in Section 2. A generic point in $E'$ is denoted by $(y, w)$. In this setting, we specialize further, assuming that

$$(4.2) \qquad E = E_1 \times E_2 \subseteq \mathbb{R} \times \mathbb{R} \quad \text{and} \quad \mathcal{E} \subseteq \mathcal{B}(\mathbb{R}^2)$$

and that $Q$ is a bivariate distribution on $\mathcal{B}(\mathbb{R}^2)$ which can be identified by a stochastic vector $(X, W)$. The generalized conditional expectation $\mu_{W|X}[g]$ is the conditional expectation of $g(X, W)$ given $X = y$, that is,

$$(4.3) \qquad \mu_{W|X}[g] \overset{\text{def}}{=} \mathbb{E}_Q[g(X, W) \mid X = y], \qquad g \in L^1(E', \mathcal{E}', Q).$$

The following definition of a generalized conditional variance is an immediate consequence of (4.3):

$$(4.4) \qquad \sigma_g^2(y) \overset{\text{def}}{=} \mu_{W|X}[g^2](y) - \mu_{W|X}^2[g](y).$$

From (4.1), it follows that $Q$ is independent of the specific normalization which identifies a particular $\pi_s$ as a function of $s$. Hence, $\mu_{W|X}$ is independent of a specific atom and also of the assumption that $m_0 = 1$, which is important in applications to real data. Suppose that $C_1'$ is an alternative to $C_1$ and let $\mu_{W|X}'$ be the alternative conditional expectation. Then

$$\mu_{W|X}[g]1_{C_1 \cap C_1'} = \mu_{W|X}'[g]1_{C_1 \cap C_1'}.$$

If $\pi_s$ is absolutely continuous with respect to two-dimensional Lebesgue measure and has density $p_s$, then with

$$p_s^X(y) \overset{\text{def}}{=} \int p_s(y, w) \ dw$$

and

$$p_{W|X}(w \mid y) \overset{\text{def}}{=} \frac{p_s(y, w)}{p_s^X(y)} 1(p_s^X(y) > 0), \qquad (y, w) \in C_1 \times E_2,$$

we have

$$\mu_{W|X}[g](y) = \int g(y, w) p_{W|X}(w \mid y) \, dw.$$

In the remainder of this section, we assume that the state space is given as in (4.2). By $x$, we denote a fixed point in $E_1$ and by $C_0$, the support of $g_W$. Let $\mathcal{M}_x(h) \stackrel{\text{def}}{=} \mathcal{N}_x(h) \times C_0$, $0 \le h \le 1$, $\mathcal{M} \stackrel{\text{def}}{=} \mathcal{M}_0(1)$, $\mathcal{M}_x \stackrel{\text{def}}{=} \mathcal{M}_x(1)$ and $\mathcal{M}\{x\} \stackrel{\text{def}}{=} \mathcal{M}_x(0)$, $\mathcal{N}_x$ being defined at the beginning of Section 2.1. If $\lambda$ is the initial measure for the compound chain, then $\lambda_{W|X}(\cdot \mid x_0)$ is the conditional distribution of $W_0$, conditional on $X_0 = x_0$, and $\lambda_X = \lambda(\cdot \times E_2)$ is the marginal initial measure for $X_0$.

### 4.2. Conditions and dependence.
In order to extend our asymptotic result to the dependent case, we will apply the conditions stated below.

The first set of conditions is related to the Markovian structure of the compound chain. Basically, we assume that the $\{X_t\}$-process also determines the $\beta$-null structure for the compound process. This holds in the independent case (cf. Lemma 3.1).

$D_1$  (i) The process $\{(X_t, W_t)\}$ is a $\phi$-irreducible, Harris recurrent Markov chain with state space given by (4.2) and transition probability function $P$.

(ii) The minorization inequality (2.2) with $(s, \nu)$ is satisfied with a corresponding invariant measure $\pi_s$.

$D_2$  (i) The marginal process $\{X_t\}$ is a $\phi_1$-irreducible, Harris recurrent Markov chain on $(E_1, \mathcal{E}_1)$ with transition probability function $P_1$.

(ii) The minorization inequality (2.2) is satisfied with $(s_1, \nu_1)$.

(iii) The Markov chain $\{X_t\}$ is $\beta$-null recurrent.

(iv) There exists a set $C_1 \in \mathcal{E}_1^+$ such that $\xi \stackrel{\text{def}}{=} 1_{C_1 \times E_2} \in \mathcal{E}^+$ is $\pi_s$-integrable.

$D_3$  (i) The invariant measure $\pi_s$ has a density, $p_s$, with respect to the two-dimensional Lebesgue measure.

(ii) $\int p_s(x, w) \, dw > 0$.

(iii) $\lim_{\delta \downarrow 0} \int |p_s(x + \delta, w) - p_s(x, w)| \, dw = 0$.

(iv) The marginal transition probability function $P_1$ is independent of any initial distribution $\lambda$.

$D_4$  (i) $g_W$ is bounded and $0 < \int |g_W(w)| \, dw < \infty$.

(ii) The set $\mathcal{N}_x \otimes C_0$ is small.

(iii) $\mu_{W|X}[g_W](x) = 0$.

$D_5$  $\forall \{A_h\} \in \mathcal{E}^\infty : \lim_{h \downarrow 0} A_h \downarrow \varnothing : \lim_{h \downarrow 0} \overline{\lim}_{y \to x} \int P((y, w), A_h) |g_W|(w) \, dw = 0$.

Conditions $D_1$–$D_3$ and $D_5$ are essentially rephrased versions of the conditions used in Theorem 3.5. Condition $D_4$ introduces stronger restrictions on $g_W$. In Section 3, the boundedness and smallness were avoided by means of a truncation technique which fit that situation. It is not obvious how to find a similar truncation procedure in the dependent situation. Possibly, the concept of *asymptotic independence* (to be introduced in Definition 4.1) could be of use. In a simulation experiment in Section 5 with an unbounded $g_W$ having noncompact support, we obtain results indicative of the asymptotics being valid under the more general conditions on $g_W$ used in Section 3.

Condition $D_4$(iii), which contains the restriction on the dependence relationship between $\{X_t\}$ and $\{W_t\}$, at first sight seems very stringent, but it will now be shown that it is, in fact, a natural extension of the type of dependence that is used in standard linear cointegration theory. Since this is important in an econometric interpretation of our results, we will consider it in some detail.

We begin by defining the concept of asymptotic independence in this context.

DEFINITION 4.1. Suppose that $\{(X_t, W_t)\}$ is a null recurrent Markov chain. The two marginal processes $\{X_t\}$ and $\{W_t\}$ are *asymptotically independent* if the invariant measure $\pi_s$ factors into a product of two measures which correspond to the $X$-component and the $W$-component.

If $\{X_t\}$ and $\{W_t\}$ are asymptotically independent, then the conditional expectation given by (4.3) reduces to a constant whenever $g(y, w) = g(w)$ and $D_4$(iii) follows if this constant is zero.

It may seem that asymptotic independence is tantamount to requiring independence, but this is not the case because having $\{X_t\}$ nonstationary (and null recurrent) and $\{W_t\}$ stationary is a special situation, where, intuitively, the "small" process $\{W_t\}$ has little influence on the "big" process $\{X_t\}$ in the long term, but allows for dependence for fixed $t$, as is the case for linear cointegration models. This phenomenon is handled more formally in the following example, which extends well-known results in linear cointegration theory (see, e.g., [15], pages 586–589).

EXAMPLE 4.1 (Asymptotic independence). In this example, we prove asymptotic independence between a random walk and a stationary autoregressive process, despite the fact that they are linked for each $t$. Moreover, we prove that conditions $D_1$–$D_5$ are satisfied. This means that the common invariant measure for these two processes factors as if the processes were independent. The processes are given by

$$
\begin{aligned}
X_t &= X_{t-1} + e_t, \\
W_t &= aW_{t-1} + be_t + u_t, \qquad |a| < 1,
\end{aligned}
\tag{4.5}
$$

where $\{e_t\}$ and $\{u_t\}$ are independent i.i.d. processes with finite third order moments and distribution functions $F_e$ and $F_u$, respectively. Moreover, we assume that these distribution functions have densities $f_e$ and $f_u$, respectively, with respect to the Lebesgue measure in $R^1$. In addition, we assume that both densities are bounded away from zero on some interval $[-c, c]$ with $c > 0$. Let $\pi_W$ denote the stationary measure for $\{W_t\}$ and $p_W$ the corresponding density. Likewise, let $\pi_X(dy) = dy$.

First, we find the density of the transition probability function for (4.5):

$$F_{X,W}(x, w \mid x_0, w_0) \stackrel{\text{def}}{=} P(X_1 \leq x, W_1 \leq w \mid X_0 = x_0, W_0 = w_0)$$

$$= P(x_0 + e_1 \leq x, aw_0 + be_1 + u_1 \leq w)$$

$$= \int P(x_0 + e_1 \leq x, aw_0 + be + u_1 \leq w \mid e_1 = e) F_e(de)$$

$$= \int 1(e \leq x - x_0) F_u(w - aw_0 - be) F_e(de)$$

and

$$f_{X,W}(x, w \mid x_0, w_0) = \frac{\partial}{\partial x} \frac{\partial}{\partial w} F_{X,W}(x, w \mid x_0, w_0)$$

$$= f_e(x - x_0) f_u(w - aw_0 - b(x - x_0)).$$

The function $f_{X,W}(x, w \mid x_0, w_0)$ is the density of the compound transition probability and from the assumption on $f_e$ and $f_u$, it follows that

$$(4.6) \qquad \inf_{(x,w,x_0,w_0) \in \mathcal{C}^4} f_{X,W}(x, w \mid x_0, w_0) > 0$$

for some Lebesgue-positive compact set $\mathcal{C}$ in $\mathbb{R}$. By (4.6), we can choose an atom $s \otimes \nu$ which is equal to a constant times $1_{\mathcal{C}} \otimes \ell I_{\mathcal{C}}$, where $\ell I_{\mathcal{C}}$ is the restriction of the Lebesgue measure to the set $\mathcal{C}$. In a similar way, we use the definitions of $\{X_t\}$ and $\{W_t\}$ to get marginal minorization inequalities, $P_i \geq s_i \otimes \nu_i$, where $P_1$ corresponds to the $X$-process and $P_2$ corresponds to the $W$-process.

If $p_s$ is an invariant density, then $p_s$ satisfies

$$(4.7) \qquad p_s(x, w) = \int p_s(x_0, w_0) f_{X,W}(x, w \mid x_0, w_0) \, dx_0 \, dw_0.$$

On the other hand, if we can find a function $p_s$ which satisfies (4.7) such that $\pi_s \stackrel{\text{def}}{=} \int p_s$ is an invariant measure with $\pi_s s = 1$, then this $p_s$ is the unique invariant density satisfying $\pi_s s = 1$. We will show that

$$(4.8) \qquad p_s(x, w) \stackrel{\text{def}}{=} c^{-1} p_{s_2}(w), \qquad c = \int s(y, w) p_{s_2}(w) \, dy \, dw$$

satisfies (4.7), where the constant $c$ is defined so that $\pi_s(s) = 1$. The measure defined by (4.8) satisfies (4.7) iff $p_s = p'_s$, where

$$(4.9) \qquad p'_s(x, w) \stackrel{\text{def}}{=} \iint c p_2(w_0) f_{X,W}(x, w \mid x_0, w_0) \, dx_0 \, dw_0.$$

From (4.9), we get

$$
\begin{aligned}
p'_s(x, w) &= \int c p_{s_2}(w_0) \left\{ \int f_{X,W}(x, w \mid x_0, w_0) \, dx_0 \right\} dw_0 \\
&= \int c p_{s_2}(w_0) \left\{ \int f_e(x - x_0) f_u(w - a w_0 - b(x - x_0)) \, dx_0 \right\} dw_0 \\
&= \int c p_{s_2}(w_0) \left\{ \int f_e(\xi) f_u(w - a w_0 - b \xi) \, d\xi \right\} dw_0 \\
&= \int c p_{s_2}(w_0) f_{W,W_0}(w \mid w_0) \, dw_0 \\
&= c p_{s_2}(w) \\
&= p_s(x, w),
\end{aligned}
$$

where we have used the fact that the transition probability density function for $\{W_t\}$ is given by

$$
f_{W,W_0}(w \mid w_0) = \int f_u(w - a w_0 - be) f_e(e) \, de.
$$

Since $p_{s_1}$ is constant, this means that $p_s(x, w) = c_1 p_{s_1}(x) p_{s_2}(w)$, where $c_1$ is a constant, hence the two marginal processes are asymptotically independent.

Let $g_W$ be any bounded, real, measurable function defined on $\mathbb{R}$ with compact support. By definition of the model, we have that D$_1$ is satisfied. Since $\{X_t\}$ is a random walk with a smooth noise process possessing a finite third order moment, the random walk is $\beta$-null recurrent. Since we have established asymptotic independence, condition D$_2$(iv) becomes trivial and, likewise, condition D$_3$. From (4.6), and since $g_W$ is assumed to be small, we infer that D$_4$ holds. The last condition, D$_5$, holds since the transition probability function is smooth. Thus, conditions D$_1$–D$_5$ are satisfied.

REMARK 4.1. The assumption on the $\{W_t\}$-process can be relaxed in this example. It is sufficient that $\{W_t\}$ is a stationary, nonlinear, autoregressive process. In (4.5), we may also replace the constant $b$ with a measurable function $\psi$, with $\mathbb{E}\psi^4(e)$ and $\sup \psi 1_{[-c,c]}$ finite. On the other hand, the calculations made in the example are based on the linearity of the $\{X_t\}$-process, with one interesting exception. Let $\{X_t\}$ be given by (4.5). Suppose that

$$
X'_t = \Phi(X_t),
$$

where $\Phi$ is a bijective measurable map between $E_1$ and $E'_1$. Then the processes $\{X'_t\}$ and $\{W_t\}$ are asymptotically independent.

Suppose that $e'_t \overset{\text{def}}{=} be_t + u_t$ is bounded. Then $\{W_t\}$ is uniformly recurrent (cf. [24], Example 5.6, page 93) and we can use the fact that $g_W(w) = w$. Imposing appropriate conditions, the uniform recurrence still holds in the nonlinear case.

The specialization in the next example makes the connection to linear cointegration even more explicit.

EXAMPLE 4.2. If the residuals in (4.5) are Gaussian, then we can calculate the conditional expectation for fixed $t$, and the rate at which we approach asymptotic independence and the fulfillment of $D_4$(iii) for Example 4.1.

$$\mathbb{E}(W_t \mid X_t) = \frac{\mathbb{E}(W_t X_t)}{\mathbb{E}(X_t^2)} X_t$$

and

$$\theta_t \stackrel{\text{def}}{=} \mathbb{E}(W_t X_t) = b\sigma_e^2 \frac{1 - a^{t+1}}{1 - a} = \frac{b\sigma_e^2}{1 - a} + \mathcal{O}(1).$$

Hence,

$$\mathbb{E}(W_t \mid X_t) = [\theta_t][t^{-1} X_t] = \mathcal{O}(1) \qquad \text{a.s.}$$

by the strong law of large numbers. Likewise, it follows that the instantaneous correlation between $X_t$ and $W_t$ decreases toward zero,

$$\text{corr}(X_t, W_t) = \mathcal{O}(t^{-1/2}).$$

However, $\{(X_t, W_t')\}$ is not Gaussian, where $W_t' = g_W(W_t)$.

4.3. *Asymptotic results.* After clarifying the relationship between various $\pi$-measures in Lemma 4.1, the main result is stated in Theorem 4.1. We denote by $\pi_{s_1}$ the invariant measure for $\{X_t\}$ implicitly defined by $D_2$ and we write $\pi_s^X$ for the $X$-marginal invariant measure of the compound chain defined by $D_3$.

LEMMA 4.1. *Assume that $D_1$ and $D_2$ are satisfied. Then the compound process is $\beta$-null recurrent and*

$$\pi_{s_1} 1_{C_1} = \frac{\pi_s 1_{C_1 \times E_2}}{\pi_s^X s_1}.$$

PROOF. Let $1_{C_1 \times E_2}$ be $\pi_s$-integrable according to $D_2$(iv). Let $C_2 \subseteq C_1$ such that $C_2$ is a small set for the $\{X_t\}$ -chain and $\xi = 1_{C_2 \times E_2} \in \mathcal{E}^+$. Since $\nu$ is a small measure and $\xi$ is $\pi_s$-integrable, the conditions in the ratio limit theorem (cf. [24], page 130) are satisfied and we get

$$\frac{\sum_{t=0}^{n} \nu P^t \xi}{\sum_{t=0}^{n} \nu P^t s} = \frac{\pi_s \xi}{\pi_s s} + \mathcal{O}(1) = \pi_s \xi + \mathcal{O}(1).$$

By $D_2(i)$, we have that

$$\nu_X P_1^t(dx) = \mathbb{P}_\nu(X_t \in dx, W_t \in E_2) = \nu P^t(dx \times E_2)$$

so that

$$\nu P^t \xi = \nu_X P_1^t 1_{C_2}.$$

Then

$$(4.10) \qquad \frac{\sum_{t=0}^n \nu_X P_1^t 1_{C_2}}{\sum_{t=0}^n \nu P^t s} = \frac{\pi_s \xi}{\pi_s s} + \mathcal{O}(1) = \pi_s \xi + \mathcal{O}(1).$$

On the other hand, since $C_2$ is small for the $\{X_t\}$-chain, we have

$$\frac{\sum_{t=0}^n \nu_X P_1^t 1_{C_2}}{\sum_{t=0}^n \nu_1 P_1^t s_1} = \frac{\pi_{s_1} 1_{C_2}}{\pi_{s_1} s_1} + \mathcal{O}(1) = \pi_{s_1} 1_{C_2} + \mathcal{O}(1).$$

Combining these two asymptotic relations gives

$$\frac{\sum_{t=0}^n \nu_1 P_1^t s_1}{\sum_{t=0}^n \nu P^t s} = \frac{\pi_s 1_{C_2 \times E_2}}{\pi_{s_1} 1_{C_2}} + \mathcal{O}(1).$$

Since the left-hand side does not depend on the actual $C_2$, it follows that

$$\pi_{s_1} 1_C = c_0^{-1} \pi_s 1_{C \times E_2}, \qquad C \in \mathcal{E}_1$$

for a fixed constant $c_0$. The constant can be expressed as $c_0 = \pi_s^X s_1$. The denominator of the fraction on the left-hand side of (4.10) has exactly the same asymptotic rate as the numerator. Hence, the compound chain is $\beta$-null recurrent since $\{X_t\}$ is $\beta$-null recurrent (cf. KT). $\square$

The following result is a modification of Theorem 3.5, which allows dependence between processes $\{W_t\}$ and $\{X_t\}$ in (3.1).

THEOREM 4.1.   *Assume* $D_1$–$D_5$. *Moreover, assume that the kernel* $K$ *satisfies* $B_1$–$B_2$ *and that for some* $\epsilon > 0$, *the bandwidth satisfies* $h_n^{-1} \ll n^{\beta - \epsilon}$. *Then for all initial measures* $\lambda$,

$$(4.11) \qquad h_n^{1/2} S_n^{1/2}(K_{x,h_n}) \left\{ \widehat{f}(x) - f(x) - \mu(g_{h_n}) \frac{S_n(s)}{S_n(K_{x,h_n})} - \frac{\pi_s^X I_{K_{x,h_n}} \psi_x}{\pi_s^X K_{x,h_n}} \right\}$$

$$\xrightarrow[n]{d} \mathcal{N}(0, \sigma_{g_W}^2(x) \|K\|_2^2),$$

*where* $\sigma_{g_W}^2(x)$ *is given by* (4.4).

   *If the density* $p_s^X$ *and the function* $f$ *possess continuous derivatives of second order, then the second bias term* $\pi_s^X I_{K_{x,h_n}} \psi_x / \pi_s^X K_{x,h_n}$ *is negligible when* $h_n^{-1} \gg n^{\beta/5+\epsilon}$. *If* $p_s^{(2,0)} = \frac{\partial^2 p_s}{dx^2}$ *exists and satisfies* $\int \overline{\lim}_{y \to x} |p_s^{(2,0)}(y,w)| \times |g_W|(w)\, dw < \infty$, *then the first bias term is negligible when* $h_n^{-1} \gg n^{\beta/5+\epsilon}$.

PROOF. The proof of this result can be seen as a modification of the proof of Theorem 3.5. That proof was built on Theorem 3.1, which, in turn, was based on C1–C7. By D1 and Lemma 4.1, the $\{(X_t, W_t)\}$ -process is $\beta$-null recurrent.

Let $g_h = K_{x,h} \otimes g_W$, $g_h^0 = \mu(g_h)s$, $\theta_h = K_{x,h} \cdot [\psi_x]$, $\theta_h^0 = K_{x,h} \cdot [\psi_x - a_h]$, $\psi_x = f - f(x)$ and $a_h = \pi_s^X I_{K_{x,h}} \psi_x / \pi_s^X K_{x,h}$. Then

$$\widehat{f}(x) = f(x) + \frac{S_n(g_h - g_h^0)}{S_n(K_{x,h})} + \frac{S_n(g_h^0)}{S_n(K_{x,h})} + \frac{S_n(\theta_h^0)}{S_n(K_{x,h})} + a_h.$$

In this notation, the left-hand side of (4.11) equals

$$(4.12) \qquad \{hS_n(K_{x,h})\}^{1/2} \left\{ \frac{S_n(g_h - g_h^0)}{S_n(K_{x,h})} + \frac{S_n(\theta_h^0)}{S_n(K_{x,h})} \right\}.$$

As noted in the proof of Theorem 3.5, it is enough to prove that

$$S_n^{-1/2}(K_{x,h_n}) h_n^{1/2} S_n(g_{h_n} - g_{h_n}^0) \xrightarrow[n]{d} \mathcal{N}(0, \sigma_{g_W}^2(x) \|K\|_2^2),$$

since the second term of (4.12) is $\mathcal{O}_P(1)$.

By D4(i)–(ii), $|g_h|$ is a small function and

$$
\begin{aligned}
(4.13) \qquad \mu(g_h) &= \pi_s g_h \\
&= \pi_s(K_{x,h} \otimes g_W) \\
&= \iint p_s(x + hu, w) K(u) g_W(w) \, dw \, du \\
&= \int p_s^X(x + hu) K(u) \mu_{W|X}[g_W](x + hu) \, du \\
&= \mathcal{O}(1),
\end{aligned}
$$

where we have used the fact that D3 implies both $p_s^X$ and $\mu_{W|X}[g_W]$ are continuous at the point $x$ and D4(iii), which ensures that the generalized conditional expectation is zero at $x$.

We also find that

$$\mu(|g_h|) = p_s^X(x) \mu_{W|X}[|g_W|](x) + \mathcal{O}(1).$$

Let $g_h' = g_h - g_h^0$. Then since $\mu(g_h') = 0$,

$$(4.14) \quad \sigma^2(g_h') = \pi_s g_h'^2 + 2h^{-1} \Delta^\star(g_h', hg_h'), \qquad \Delta^\star(g_h, f_h) \overset{\text{def}}{=} \pi_s I_{g_h} P G_{s,\nu} f_h,$$

using (A.12). By (4.13) and (4.14),

$$h\pi_s g_h'^2 + 2\Delta^\star(g_h', hg_h') = h\pi_s g_h^2 + 2\Delta^\star(g_h, hg_h) + \mathcal{O}(1).$$

Hence, an asymptotic variance, $\sigma^2 \overset{\text{def}}{=} \lim_{h \downarrow 0} h\sigma^2(g_h')$, if it exists, is given by

$$(4.15) \qquad \sigma^2 = \lim_{h \downarrow 0} \{ h\pi_s g_h^2 + 2\Delta^\star(g_h, hg_h) \}.$$

In order to verify (4.15), we begin by showing that the first term on the right-hand side of (4.15) satisfies

$$h\pi_s(g_h^2) = \|K\|_2^2 \sigma_{gW}^2(x) p_s^X(x) + \mathcal{O}(1),$$

where the conditional variance is given by (4.4).

Indeed, by the definition of $g_h$, we find that

$$\pi_s(g_h^2) = \pi_s(K_{x,h}^2 \otimes g_W^2)$$

$$= \int p_s(y,w) K_{x,h}^2(y) g_W^2(w) \, dy \, dw$$

$$= h^{-1} \int p_s(x+hu,w) K^2(u) g_W^2(w) \, du \, dw$$

$$= h^{-1} \left[ \|K\|_2^2 \int p_s(x,w) g_W^2(w) \, dw + \mathcal{O}(1) \right]$$

$$= h^{-1} p_s^X(x) [\|K\|_2^2 \sigma_{gW}^2(x) + \mathcal{O}(1)].$$

The next task is to show that $\Delta^\star(g_h, hg_h)$ is asymptotically negligible. Let

$$f_h(y,w) \stackrel{\text{def}}{=} [hK_{x,h}(y) - 1_{\{x\}}(y) K(0)] g_W(w)$$

$$= [1_{\{x\}^c} \cdot hK_{x,h} \otimes g_W](y,w)$$

so that

$$hg_h = f_h + K(0)[1_{\{x\}} \otimes g_W] = f_h + \delta_h,$$

say.

By D$_3$(i), we find that $\pi_s|\delta_h| = 0$ and thus $\pi_s PG_{s,\nu}\delta_h = 0$. Hence,

$$\Delta^\star(g_h, hg_h) = \Delta^\star(g_h, f_h)$$

$$= \iint p_s(x+hu,w) K(u) g_W(w) PG_{s,\nu} f_h(x+hu,w) \, du \, dw,$$

where we have inserted the invariant density and made a standard substitution.

Let

$$\eta_h = \|K\|_\infty G_{s,\nu} \{1_{\{x\}^c} \cdot \mathcal{N}_x(h)\} \otimes |g_W|$$

so that

(4.16)
$$\begin{aligned} &|\Delta^\star(g_h, f_h)| \\ &\leq \iint p_s(x+hu,w) K(u) |g_W(w)| P\eta_h(x+hu,w) \, du \, dw. \end{aligned}$$

By D$_4$(i)–(ii) and Nummelin ([24], Proposition 5.13, page 80), the function $\eta_h$ is bounded. Since $\{1_{\{x\}^c} \cdot \mathcal{N}_x(h)\} \otimes |g_W| \downarrow 0$ pointwise, $\lim_{h\downarrow 0} \eta_h(y,w) \downarrow 0$.

Let $\epsilon > 0$ and $A_h \stackrel{\text{def}}{=} \{\eta_h > \epsilon\}$. Then $\{A_h\}$ satisfies $D_5$. Inserting $\eta_h = I_{A_h^c}\eta_h + I_{A_h}\eta_h$ into (4.16), we get

$$|\Delta^\star(g_h, f_h)| \leq \epsilon \iint p_s(x+hu,w)K(u)|g_W(w)|\,du\,dw$$

$$+ \|\eta_h\|_\infty \iint p_s(x+hu,w)K(u)|g_W(w)|$$

$$\times P((x+hu,w), A_h)\,du\,dw.$$

The main part of the last term of the above expression is bounded by

$$(4.17) \quad \int p_s(x+hu,w)K(u)\left\{\sup_{|x-y|<\epsilon_h}\int P((y,w),A_h)|g_W(w)|\,dw\right\}du$$

for all $\epsilon_h = h\sup\{|u| : u \in \mathcal{N}_0\}$ and $\lim_{h\downarrow 0}\epsilon_h \downarrow 0$.

Using $D_4$, it follows that (4.17) is $\mathcal{O}(1)$ with respect to $h$. Putting all of this together, it is clear that

$$\lim_{\epsilon\downarrow 0}\overline{\lim_{h\downarrow 0}}|\Delta^\star(g_h, f_h)| = 0.$$

Thus, we have so far proved that

$$h\sigma^2(g_h) = h\pi_s(g_h^2) + \mathcal{O}(1) = \|K\|_2^2\sigma_{g_W}^2(x)p_s^X(x) + \mathcal{O}(1).$$

We must also check $h\sigma^2(|g_h'|)$. This quantity is given by

$$h\sigma^2(|g_h'|) = h\pi_s|g_h'|^2 + 2\Delta^\star(|g_h'|, h|g_h'|) - h\pi_s^2|g_h'| - 2h\pi_s|g_h'|\,\pi_s(s\cdot|g_h'|),$$

by (A.12) of Appendix A. Since $\pi_s|g_h'| \leq \pi_s|g_h| + |\mu_{g_h}| = \pi_s|g_h| + \mathcal{O}(1)$ and $\pi_s|g_h| = \mathcal{O}(1)$, we have

$$h\sigma^2(|g_h'|) = h\pi_s|g_h'|^2 + 2\Delta^\star(|g_h|, h|g_h|) + \mathcal{O}(1).$$

By the same arguments as those given above, we find that

$$h\sigma^2(|g_h'|) = h\pi_s|g_h'|^2 + \mathcal{O}(1)$$
$$= h\pi_s|g_h|^2 + \mathcal{O}(1).$$

Since $g_h$ is small, we easily find that (cf. Theorem A.1 in Appendix A)

$$\mathbb{E}\|U(g_h) - \mu(g_h)\|^{2m} \leq d_m h^{-2m+1}, \qquad m \geq 1$$

and

$$\mathbb{E}\|U(|g_h| - \mu(|g_h|)\|^{2m} \leq d_m' h^{-2m+1}, \qquad m \geq 1.$$

Moreover, we have $h|g_h| \leq g_0$, $g_0 \stackrel{\text{def}}{=} c_0 1_{\mathcal{M}_x}$ and $\mathbb{P}_\lambda(U_0(g_0) < \infty) = 1$.

Thus, the assumptions in Theorem 3.1 are satisfied and (4.11) holds. It is straightforward to verify that the bias terms are negligible under the given conditions (cf. KT).  □

REMARK 4.2. If $\int p_s(y, w)g_W(w)\,dw \equiv 0$, then $\mu(g_h) \equiv 0$, by D$_3$(iii). If this assumption holds, then the stochastic bias correcting term in (4.11) is zero. If D$_4$(iii) is strengthened so as to also require asymptotic independence, then $\sigma_W(x) \equiv \sigma_W^2 = \mathbb{E}g_W^2(W_t)$.

## 5. Simulations and finite sample behavior.

Estimates similar to that in (1.1) have appeared in the cointegration literature. Our contribution, which we believe to be new, is that we have singled out classes of processes and assumptions for which an asymptotic theory of these estimates can be constructed, such that it should be possible to work out confidence intervals and bands (and possibly rigorous tests of nonlinear cointegration, in the sense discussed in this paper).

The purpose of this section is to illustrate the small-sample properties of the estimator $\hat{f}(x)$ defined by (1.1), using simulations,

A problem not encountered in the stationary case is that the simulated realizations may cover very different $x$-regions. Hence, for a fixed $x = x'$, close to the starting value $X_0 = 0$, say, of each realization, some realizations may have many observations in the neighborhood of $x'$, whereas other realizations may have none in the vicinity of $x'$ for the sample size we are considering. This kind of behavior does not occur in the stationary case, where the expected time until the process reaches $x'$ is always finite and, in practice, small when $|x'|$ is small. This means that in a finite-sample approximation of the asymptotics, we can either keep $x$ fixed and wait until we have sufficiently many observations close to $x$ or we can choose a central realization-dependent value of $x$ (e.g., the modal value of the sample) for studying the normalized ratio (3.2) of Theorems 3.5 and 4.1. We have chosen to adopt both procedures, although, clearly, we introduce some extraneous stochastics into the problem in the latter case.

A difficult and largely unresolved problem is that of choosing a proper bandwidth. Theorem 3.5 and Theorem 4.4 of KT only give the allowable rate as $n$ tends to infinity. It should be noted that these rates are different from those in the stationary case, $n$ effectively being replaced by $n^\beta$. In practice, we have found it useful to use cross-validation and to let the bandwidth $h$ depend on $x$. In fact, we have typically let $h_n$ be proportional to $\{T_C(n)\hat{p}_C(x)\}^{-1/5}$, where $\hat{p}_C(x)$ could be thought of as the locally estimated density and where it is known from KT (Lemma 3.4) that $T_C(n)$ essentially behaves as $n^\beta$.

The approximation to normality as a function of sample size, for the quantity

$$(5.1) \qquad \left[\frac{h_n \sum K_{x,h_n}}{\int K^2(u)\,du}\right]^{1/2}[\hat{f}(x) - x]$$

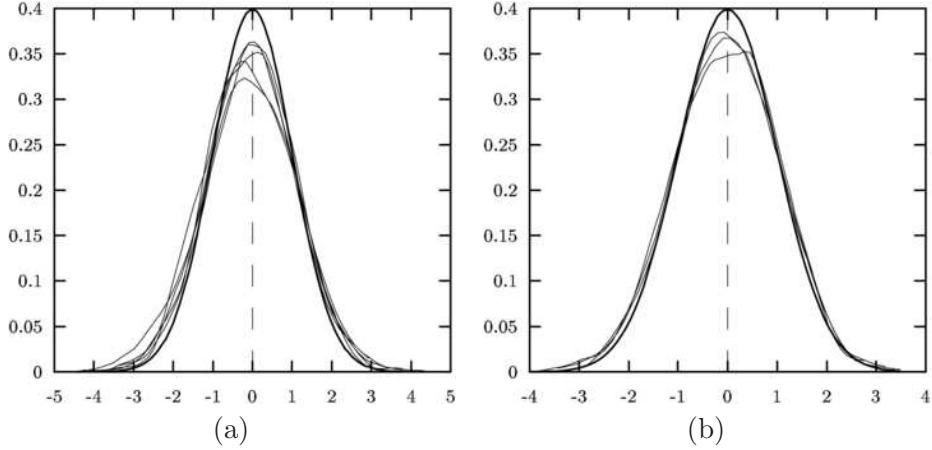

Fig. 1. (a) Thick line: *The standard normal pdf.* Thin lines: *The estimated pdfs for the quantity* $(h_n \sum K_{x,h_n} / \int K^2(u)\,du)^{1/2} [\widehat{f}(x) - f(x)]$, *at the point* $x = 7.5$, *derived from the cointegrated system* $(X_t, Z_t; t \geq 1)$, *where* $X_t = X_{t-1} + e_t$, $Z_t = f(X_t) + \varepsilon_t$, $e_t$ *and* $\varepsilon_t$ *are independent i.i.d.* $\mathcal{N}(0,1)$ *variables and* $f(x) = x$ *for all real* $x$. *The quantity is estimated by 1000 realizations and a particular realization is admitted into the evaluation as, respectively,* 100, 200, 300, 500 *and* 800 *observations are accumulated in the interval* $(5, 10)$. (b) Thick line: *The standard normal pdf.* Thin lines: *The estimated pdfs for the same quantity as in* (a), *but where a particular realization is admitted into the evaluation at the modal value. The length of the time series is* 500, 1000 *and* 3000, *respectively.*

derived from the simple cointegrated system

$$X_t = X_{t-1} + e_t, \qquad Z_t = X_t + W_t, \qquad e_t \text{ and } W_t \text{ independent} \sim \mathcal{N}(0,1)$$

at the point $x = 7.5$, is shown in Figure 1(a). 1000 realizations have been used and a particular realization is admitted into the evaluation as, respectively, 100, 200, 300, 500 and 800 observations are accumulated in the interval $(5, 10)$. For Figure 1(b), on the other hand, a fixed point $x$ has not been used; rather, $x$ has been taken to be the modal value and is thus varying from one realization to another. In this case, the length of the time series is 500, 1000 and 3000, respectively.

In Figures 2(a) and 2(b), we have considered (5.1) for the system

$$(5.2) \quad \begin{aligned} X_t &= X_{t-1} + e_t, \qquad Z_t = X_t + W_t, \qquad W_t = \sqrt{0.5}e_t + \sqrt{0.5}\varepsilon_t, \\ \varepsilon_t &\text{ and } e_t \text{ independent} \sim \mathcal{N}(0,1) \end{aligned}$$

to test the asymptotics in the case of dependence between $\{X_t\}$ and $\{W_t\}$, as described in Section 4. As in Figure 1(b), $x$ is taken to be the modal value in Figure 2(b). For both Figure 1 and Figure 2, it is seen that the finite sample distribution gets reasonably close to the asymptotic normal distribution.

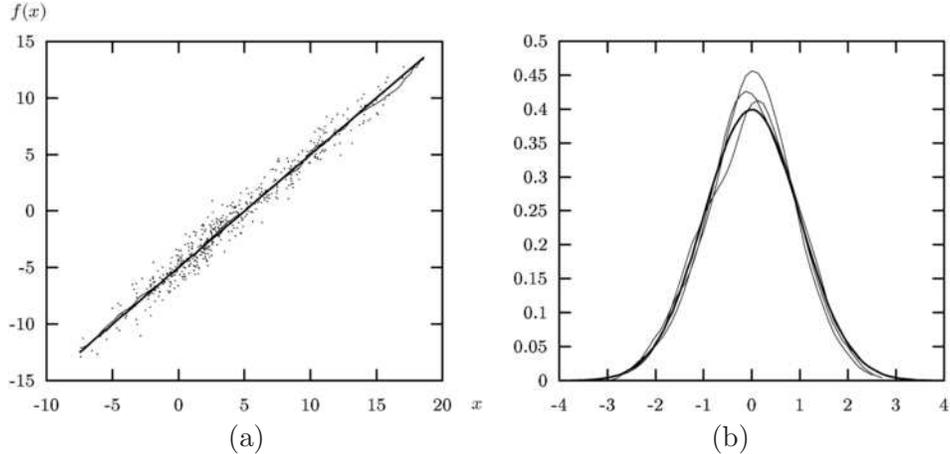

(a)                                          (b)

Fig. 2.  (a) Thick line: *The true transfer function* $f(x) = x - 5$. *Dots are* $z_t$ *plotted against* $x_t$, $t \geq 1$. *We have 500 observations. Thin line: Estimated transfer function* $\hat{f}$, *built on 500 observations from the cointegrated system* $(X_t, Z_t; t \geq 1)$, *where* $X_t = X_{t-1} + e_t$ *and* $Z_t = f(X_t) + W_t$, *where* $W_t = \sqrt{0.5}e_t + \sqrt{0.5}\varepsilon_t$ *for* $t \geq 1$. $(e_t, \varepsilon_t)$ *are i.i.d.* $\mathcal{N}(0, I)$ *vectors for* $t \geq 1$ *and* $I$ *is the identity matrix. Finally,* $f(x) = x - 5$ *for* $x$ *real. (b) Thick line: The standard normal pdf. Thin lines: The estimated pdfs for the quantity* $(h_n \sum K_{x,h_n}/\int K^2(u)\,du)^{1/2}[\hat{f}(x) - f(x)]$, *derived from the cointegrated system* $(X_t, Z_t; t \geq 1)$, *where* $X_t - X_{t-1} + e_t$ *and* $Z_t = f(X_t) + W_t$, *where* $W_t = \sqrt{0.5}e_t + \sqrt{0.5}\varepsilon_t$ *for* $t \geq 1$. $(e_t, \varepsilon_t)$ *are i.i.d.* $\mathcal{N}(0, I)$ *vectors for* $t \geq 1$ *and* $I$ *is the identity matrix. Finally,* $f(x) = x - 5$ *for* $x$ *real. The quantity is estimated by 1000 realizations and a particular realization is admitted into the evaluation at the modal value. The length of the time series is* 500, 1000 *and* 3000, *respectively.*

Note that $\{X_t\}$ and $\{W_t\}$ in (5.2) are asymptotically independent with $\mu_{W|X}[g_W](x) = 0$ and $\sigma_{g_W}(x) = 1$. Actually, in (5.2), $g_W(W_t) = W_t$ and the assumptions D₄(i) and D₄(ii) are not satisfied, this being something we wanted to test by means of this simulation experiment. On the other hand, with the exception of the independence assumption, the other assumptions in Theorem 3.5 are satisfied. We also carried out an experiment with

$$W_t = \tfrac{1}{\sqrt{3}}e_t + \tfrac{1}{\sqrt{3}}e_{t-1} + \tfrac{1}{\sqrt{3}}\varepsilon_{t-1}.$$

In this case, $\{(X_t, W_t)\}$ is not Markov. On the other hand, $\{(X_t, W_t, e_t)\}$ is Markov and $\{X_t\}$ is asymptotically independent of the Markov process $\{(W_t, e_t)\}$. The distributional results were similar to those of Figure 2(b). More simulation experiments and a real data example are given in [20].

**6. Some final remarks on nonlinear cointegration.** This paper can be looked at in two ways: (i) it is an attempt to establish a statistical theory for nonparametric regression with a nonstationary regressor and (ii) in addition, it is seeking to relate this framework to the problem of nonlinear cointegration. There are a host of open problems for both. For (ii), it is of particular interest to weaken conditions D₄(ii) and D₄(iii) on $g_W$, alternatively, letting $W'_t = g_W(X_t, \ldots, X_{t-p}, W_t)$, as indicated in the second paragraph of Section 3. But, there are also conceptual issues involved concerning the function $f$. In a nonparametric approach like ours, $f$ is determined by the data and if $\{Z_t\}$ and $\{X_t\}$ are close to being linearly cointegrated, one expects the nonparametric estimate $\widehat{f}$ to be close to a linear function and might think that the difference between $\widehat{f}$ and a linear function could be used to test for linearity of the cointegration. One could also test for appropriate parametric functions for $f$. For the estimation of nonlinear parametric regression functions using local time arguments, see [25]. Clearly, not every parametric $f$ makes sense from a cointegration framework. As an extreme case, consider $f(x) = $ constant. Then $\{Z_t\}$ will be stationary and unrelated to $\{X_t\}$. In a cointegration framework, $\{Z_t\}$ should be nonstationary and the question arises as to whether it is possible to construct nontrivial functions $f$ such that $\{Z_t\}$ is stationary, even though $\{X_t\}$ is nonstationary. Another question is whether all such functions $f$ (e.g., the sine function) will be economically meaningful.

One of the referees has pointed out that the function $f$ may include a constant term. But, a deterministic term depending on the time parameter (e.g., a linear trend) is not included in the model. An extension of the model in this direction introduces challenging problems concerning properties of estimates for both $f$ and the trend. It seems to be quite clear that additional assumptions on $\{X_t\}$ are required since null recurrence itself is not related to the growth rate of a linear trend. The situation is much more specific in the random walk case, where the variance of $\{X_t\}$ increases linearly and $f$ is linear.

Still another issue is whether $f$ should be required to be one-to-one for it to be meaningful in a cointegration framework. Requiring $f$ to be one-to-one has the advantage of allowing the possibility of expressing $\{X_t\}$ in terms of $\{Z_t\}$, making for a more symmetric relationship. To estimate such an inverse relationship would be nontrivial since it would require an extension of the theory to the case where the regressor is a function of a Markov chain.

In the linear cointegration case, the concept of cointegration is intimately connected with the so-called error correction representation (cf. [18]). Nonlinear extensions have centered on both nonlinear error correction and nonlinear cointegration (see, e.g., [7, 8, 13, 17]). It remains to explore possible connections between these models and the approach presented in this paper.

Nonlinear cointegration extensions are more demanding and are at the core of the present paper. Only a few attempts of such an extension can be found in the literature. Specific nonlinear cointegration relationships in terms of threshold models have been studied by Hansen and Seo [16] and Bec and Rahbek [3]. Escribano and Mira [8] suggest definitions of $I(0)$ and $I(1)$ which are useful in a parametric nonlinear context and study several large- and small-sample properties of nonlinear least squares estimation. Related work in nonlinear parametric regression theory has appeared in Park and Phillips [25]. Nonparametric estimates of nonlinear cointegration have been computed from data by Granger and Hallman [12] and Aparicio and Escribano [1]. However, no attempt has been made to study the asymptotic properties of nonparametric estimators either for nonlinear error correction or cointegration models.

## APPENDIX A

In this appendix, we assume that $\{X_t\}$ is an aperiodic, $\phi$-irreducible Markov chain with state space $(E, \mathcal{E})$, where $\mathcal{E}$ is countably generated. We also assume that the transition probability $P$ satisfies the minorization inequality, (2.2), that is, $P \geq s \otimes \nu$, and that $\{X_t\}$ is Harris recurrent. Recall the taboo transition probability $H = P - s \otimes \nu$, the taboo kernel $G_{s,\nu} = \sum_{j=0}^{\infty} H^j$, the index set $\Delta_r^m = \{\alpha \in \mathcal{N}_+^r : \sum_{j=1}^r \alpha_j = m\}$, the multinomial coefficient $\binom{m}{\alpha} = \frac{m!}{\alpha_1! \cdots \alpha_r!}$ and the moment function

$$(A.1) \qquad \psi_{r,\alpha} = \sum_{j \in \mathcal{N}_{0,+}^r} H^{j_1} I_{g_{j_1}^{\alpha_1}} \cdots H^{j_r} I_{g_{j_1+\cdots+j_r}^{\alpha_r}} 1.$$

The $r$-Cartesian product of the set of integers where all but the first coordinate are strictly positive is denoted by $\mathcal{N}_{0,+}^r$.

**A.1. Higher-order moments.** An expression for the moments of a $U$-block is derived from a moment formula for a real sequence (cf. [21]):

LEMMA A.1. *Let $\{a_t\}$ be a real sequence and $m \geq 1$ an integer. Then*

$$\left\{\sum_{k=0}^n a_k\right\}^m = \sum_{r=1}^m \sum_{\alpha \in \Delta_r^m} \binom{m}{\alpha} \sum_{j \in \mathcal{N}_{0,+}^r} 1(s_r \leq n) a_{s_1}^{\alpha_1} \cdots a_{s_r}^{\alpha_r},$$

*where $s_\ell = j_1 + \cdots + j_\ell$, $\ell = 1, \ldots, r$.*

THEOREM A.1. *Let $\underline{g} = \{g_t\}$ be a sequence of real-valued measurable functions defined on $E$. Let $U_0 = U_0(\underline{g}) \stackrel{\text{def}}{=} \sum_{k=0}^{\tau} g_k(X_k)$. Then*

$$\mathbb{E}_x U_0^m = \sum_{r=1}^m \sum_{\alpha \in \Delta_r^m} \binom{m}{\alpha} \psi_{r,\alpha}(x).$$

REMARK A.1. If $g_k \equiv g$, then

$$\psi_{r,\alpha} = \sum_{j \in \mathcal{N}^r} H^{j_1} I_{g^{\alpha_1}} \cdots H^{j_r+1} I_{g^{\alpha_r}} 1 = [G_{s,\nu} I_{g^{\alpha_1}} H] \cdots [G_{s,\nu} I_{g^{\alpha_{r-1}}} H] G_{s,\nu} I_{g^{\alpha_r}} 1.$$

PROOF OF THEOREM A.1. The main ingredient in this proof is the lemma formulated above, together with the Markov property.

Let $\mathcal{B}_s = \mathcal{F}_s^X \vee \mathcal{F}_{s-1}^Y$, $A_s = (1 - Y_s)$, $B_{0,s} = \{\tau \geq s\} = \prod_{k=0}^{s-1} A_k$ and $B_{t,t+h} = \prod_{k=t}^{t+h-1} A_k$. From Lemma A.1, the definition of $\Delta_r^m$ and $a_k = g(X_k)1(\tau \geq k)$ with $n = \infty$, we get

$$U_0^m = \sum_{r=1}^{m} \sum_{\alpha \in \Delta_r^m} \binom{m}{\alpha} J_{r,\alpha}, \qquad J_{r,\alpha} = \sum_{j \in \mathcal{N}_{0,+}^r} Z_{j,\alpha},$$

with

(A.2) $$Z_{j,\alpha} \stackrel{\text{def}}{=} g_{s_1}^{\alpha_1}(X_{s_1}) \cdots g_{s_r}^{\alpha_r}(X_{s_r}) B_{0,s_r}.$$

Let $r$ and $\alpha$ be fixed and $f_k(x) \stackrel{\text{def}}{=} g_k^{\alpha_k}(x)$ for $k = 1, \ldots, r$. Then

$$J_r = J_{r,\alpha} = \sum_{j \in \mathcal{N}_{0,+}^r} Z_j, \qquad Z_j = f_{s_1}(X_{s_1}) \cdots f_{s_r}(X_{s_r}) B_{0,s_r}.$$

It is enough to prove that

(A.3) $$\mathbb{E}_x J_r = \sum_{j \in \mathcal{N}_{0,+}^r} H^{j_1} I_{f_{j_1}} \cdots H^{j_r} I_{f_{j_1 + \cdots + j_r}} 1(x)$$

for arbitrary $r$ and $\{f_k\}$. We will prove this by induction on $r$. When $r = 1$,

$$J_1 = \sum_{j_1=0}^{\infty} f_{j_1}(X_{j_1}) B_{0,j_1} = \sum_{j_1=0}^{\infty} f_{j_1}(X_{j_1}) 1(\tau \geq j_1)$$

and

$$\mathbb{E}_x J_1 = \sum_{j_1=0}^{\infty} \mathbb{E}_x f_{j_1}(X_{j_1}) 1(\tau \geq j_1) = \sum_{j_1=0}^{\infty} H^{j_1} I_{f_{j_1}} 1(x),$$

which shows that (A.3) is true for $r = 1$.

Assume that (A.3) is true for $r - 1$. Corresponding to the induction hypothesis, let $\hat{j} = (j_1, \ldots, j_{r-1})$ and $\underline{s} = s_{r-1}$. Then

(A.4) $$Z_j = f_{s_1}(X_{s_1}) \cdots f_{s_r}(X_{s_r}) B_{0,s_r} = Z_{\hat{j}} A_{\underline{s}} f_{\underline{s}+j_r}(X_{\underline{s}+j_r}) B_{\underline{s}+1, \underline{s}+j_r},$$

where $Z_{\widehat{\jmath}}A_{\underline{s}}$ is measurable $\mathcal{B}_{\underline{s}+1}$. Taking conditional expectation with respect to $\mathcal{B}_{\underline{s}+1}$ in the last part of (A.4) gives

$$
\begin{aligned}
\mathbb{E}\{f_{\underline{s}+j_r}(X_{\underline{s}+j_r})B_{\underline{s}+1,\underline{s}+j_r} \mid \mathcal{B}_{\underline{s}+1}\} &= \mathbb{E}\{f_{\underline{s}+j_r}(X_{\underline{s}+j_r})B_{\underline{s}+1,\underline{s}+j_r} \mid \mathcal{B}_{\underline{s}+1}\} \\
&= \mathbb{E}_{X_{\underline{s}+1}}\{f_{\underline{s}+j_r}(X_{j_r-1})B_{0,j_r-1}\} \\
&= H^{j_r-1}f_{\underline{s}+j_r}(X_{\underline{s}+1}),
\end{aligned}
$$

so that

$$
\begin{aligned}
\phi^0_{\underline{s}}(X_{\underline{s}+1}) &\stackrel{\text{def}}{=} \sum_{j_r=1}^{\infty} \mathbb{E}\{f_{\underline{s}+j_r}(X_{\underline{s}+j_r})B_{\underline{s}+1,\underline{s}+j_r} \mid \mathcal{B}_{\underline{s}+1}\} \\
&= \left\{ \sum_{j_r=1}^{\infty} H^{j_r-1}f_{\underline{s}+j_r} \right\}(X_{\underline{s}+1}).
\end{aligned}
\tag{A.5}
$$

Combining (A.4)–(A.5), we obtain

$$
\begin{aligned}
\mathbb{E}_x J_r &= \mathbb{E}_x \sum_{\widehat{\jmath} \in \mathcal{N}_{0,+}^{r-1}} Z_{\widehat{\jmath}} A_{\underline{s}} \phi^0_{\underline{s}}(X_{\underline{s}+1}) \\
&= \mathbb{E}_x \sum_{\widehat{\jmath} \in \mathcal{N}_{0,+}^{r-1}} Z_{\widehat{\jmath}} \mathbb{E}[A_{\underline{s}} \phi^0_{\underline{s}}(X_{\underline{s}+1}) \mid \mathcal{B}_{\underline{s}}] \\
&= \mathbb{E}_x \sum_{\widehat{\jmath} \in \mathcal{N}_{0,+}^{r-1}} Z_{\widehat{\jmath}} \phi_{\underline{s}}(X_{\underline{s}}),
\end{aligned}
\tag{A.6}
$$

where

$$
\phi_{\underline{s}}(X_{\underline{s}}) \stackrel{\text{def}}{=} \mathbb{E}\{A_{\underline{s}} \phi^0_{\underline{s}}(X_{\underline{s}+1}) \mid \mathcal{B}_{\underline{s}}\} = \mathbb{E}_{X_{\underline{s}}}\{\phi^0_{\underline{s}}(X_1)(1-Y_0)\} = H\phi^0_{\underline{s}}(X_{\underline{s}}).
$$

The conditional step above reduces the dimension of $j$ and it remains to verify that (A.3) is correct when we apply the induction hypothesis. We look at $f_{\underline{s}}(x) = f_{s_{r-1}}$ defined in (A.4). Let $f^0_{\underline{s}}(x) = f_{\underline{s}}(x)H\phi^0_{\underline{s}}(x)$. Then

$$
I_{f^0_{\underline{s}}}1(x) = f^0_{\underline{s}}(x) = I_{f_{\underline{s}}}\left[\sum_{j_r=1}^{\infty} H^{j_r}f_{\underline{s}+j_r}\right](x).
\tag{A.7}
$$

By (A.6), the product (A.2) is reduced from $r$ to $r-1$ since, using the fact that $\underline{s} = s_{r-1}$, we have

$$
\mathbb{E}_x J_r = \mathbb{E}_x \sum_{\widehat{\jmath} \in \mathcal{N}_{0,+}^{r-1}} f_{s_1}(X_{s_1}) \cdots f_{s_{r-2}}(X_{s_{r-2}}) f^0_{s_{r-1}}(X_{s_{r-1}}).
\tag{A.8}
$$

Hence, by (A.8), we can evaluate the expectation of $J_r$ by the induction hypothesis, which, together with (A.7), gives (A.3).  $\square$

COROLLARY A.1. *Let $U_0(\underline{a}, g) = \sum_{k=0}^{\tau} a_k g(X_k)$. Then*

$$\mathbb{E}_x U_0^m(\underline{a}, g) = \sum_{r=1}^{m} \sum_{\alpha \in \Delta_r^m} \binom{m}{\alpha} \psi_{r,\alpha}(x),$$

(A.9)

$$\psi_{r,\alpha} = \sum_{j \in \mathcal{N}_{0,+}^r} d_{j,\alpha}(x) \{ a_{j_1}^{\alpha_1} a_{j_1+j_2}^{\alpha_2} \cdots a_{j_1+\cdots+j_r}^{\alpha_r} \},$$

*where $d_{j,\alpha} = H^{j_1} I_{g^{\alpha_1}} \cdots H^{j_r} I_{g^{\alpha_r}} 1$. In particular, for $m = 1, 2$, we have that*

$$\mathbb{E}_x U_0(\underline{a}, g) = \sum_{j=0}^{\infty} d_j(x) a_j, \qquad d_j = H^j I_g 1,$$

(A.10)

$$\mathbb{E}_x U_0^2(\underline{a}, g) = \sum_{j=0}^{\infty} d_{j,0}(x) a_j^2 + 2 \sum_{j=0}^{\infty} \sum_{\ell=1}^{\infty} d_{j,\ell}(x) \{ a_j a_{j+\ell} \},$$

$$d_{j,\ell} = H^j I_g H^\ell I_g 1.$$

PROOF. We obtain (A.9) from (A.1) with $g_j = a_j g$. With $m = 2$, we have

$$\mathbb{E}_x U_0^2(\underline{a}, g) = \sum_{\alpha \in \Delta_1^2} \binom{2}{\alpha} \psi_{r,\alpha}(x) + \sum_{\alpha \in \Delta_2^2} \binom{2}{\alpha} \psi_{r,\alpha}(x)$$

$$= \psi_{1,2}(x) + 2\psi_{2,(1,1)}(x)$$

$$= \left[ \sum_{j_1=0}^{\infty} H^{j_1} I_{g_{j_1}^2} \right] 1(x) + 2 \left[ \sum_{j_1=0}^{\infty} H^{j_1} I_{g_{j_1}} \right] \left[ \sum_{j_2=1}^{\infty} H^{j_2} I_{g_{j_1+j_2}} \right] 1(x)$$

$$= \sum_{j=0}^{\infty} H^j g_j^2(x) + 2 \sum_{j=0}^{\infty} \sum_{s=1}^{\infty} H^j I_{g_j} H^s g_{j+s}(x).$$

Hence,

$$\mathbb{E}_\nu U_0^2(\underline{a}, g) = \sum_{j=0}^{\infty} a_j^2 \{ \nu H^j g^2 \} + 2 \sum_{j=0}^{\infty} \sum_{s=1}^{\infty} a_j a_{j+s} \{ \nu H^j I_g H^s g \}. \qquad \square$$

REMARK A.2. In particular, if $a_j \equiv 1$, we write

$$U_0 = U_0(g) = \sum_{k=0}^{\tau} g(X_k)$$

and (A.10) gives the formulas $\mathbb{E}_\nu U_0(g) = \pi_s g$ and

(A.11)

$$\mathbb{E}_\nu U_0^2(g) = \pi_s g^2 + 2\pi_s I_g H G_{s,\nu} g$$

$$= \pi_s g^2 + 2\pi_s I_g P G_{s,\nu} g - 2\pi_s I_s g \pi_s g.$$

Let $\mu(g) = \mathbb{E}_\nu U_0(g)$ and $\sigma^2(g) = \text{Var}(U(g))$. Then

$$(A.12) \quad \mu(g) = \pi_s g, \qquad \sigma^2(g) = \pi_s g^2 - \pi_s^2 g + 2\pi_s I_g P G_{s,\nu} g - 2\pi_s I_s g \pi_s g.$$

### A.2. Moment inequality.

LEMMA A.2. *Assume that* (2.2) *holds. Let* $p > 1$ *and* $\eta \in (0,1)$ *and let* $f$ *be a real-valued measurable function defined on* $E$. *Then for any probability measure* $\lambda$,

$$\lambda G_{s,\nu}^t |f| \leq c_2 \mathbb{E}_\lambda^{t/p} \{\tau^{1+2(p-1)}\} \sup_{j \geq 0} \lambda^{t/q} P^j |f|^q,$$

$$t = p/(1+\eta(p-1)), \ q = p/\eta(p-1),$$

*where* $c_2$ *is a universal constant dependent only on* $p$ *and* $\eta$.

PROOF. Let $q' = p/(p-1)$, $r = q'/(1-\eta)$, $w = 1/q'$, $v = 1/q'$ and $u = 2/q'$. Then $u = v + w$, $p^{-1} + q^{-1} + r^{-1} = 1$, $1/t = 1/p + 1/q$, $pu = 2(p-1)$ and $qv = \eta^{-1}$.

From the right-hand side of (2.10) and by the Hölder inequality, we obtain

$$G_{s,\nu}^t |f|(x) = \left[ \sum_j \mathbb{E}_x \{1(\tau \geq j) |f|(X_j)\} \right]^t$$

$$\leq \left[ \sum_j \mathbb{P}_x^{1/p}(\tau \geq j) \mathbb{E}_x^{1/q} \{|f|^q (X_j)\} \right]^t$$

$$\leq \left[ \sum_j \mathbb{P}_x^{1/p}(\tau \geq j) \{P^j |f|^q (x)\}^{1/q} \right]^t$$

$$= \left[ \sum_j (j^u \mathbb{P}_x^{1/p}(\tau \geq j)) (j^{-v} \{P^j |f|^q (x)\}^{1/q}) (j^{-w}) \right]^t$$

$$\leq \left( \sum_j j^{up} \mathbb{P}_x(\tau \geq j) \right)^{t/p} \left( \sum_j j^{-vq} P^j |f|^q (x) \right)^{t/q} \left( \sum_j j^{-wr} \right)^{t/r}$$

$$= c_1 V^{t/p} Z^{t/q},$$

say. We apply the Hölder inequality again with $p_1 = p/t$ and $q_1 = q/t$. This gives

$$\lambda G_{s,\nu}^t |f| \leq c_1 \lambda V^{t/p} Z^{t/q}$$

$$\leq c_1 [\lambda^{t/p} V] [\lambda^{t/q} Z]$$

$$= c_1 \left[ \sum_{j=0}^{\infty} j^{up} \mathbb{P}_\lambda(\tau \geq j) \right]^{t/p} \left[ \sum_{j=0}^{\infty} j^{-vq} \lambda P^j |f|^q \right]^{t/q}$$

$$= c_1 \left[ \sum_{j=0}^{\infty} j^{2(p-1)} \mathbb{P}_\lambda(\tau \geq j) \right]^{t/p} \left[ \sum_{j=0}^{\infty} j^{-\eta-1} \lambda P^j |f|^q \right]^{t/q}$$

$$\leq c_2 \left[ \sum_{j=0}^{\infty} j^{2(p-1)} \mathbb{P}_\lambda(\tau \geq j) \right]^{t/p} \sup_{j \geq 0} \lambda^{t/q} P^j |f|^q$$

$$\leq c_2 \mathbb{E}_\lambda^{t/p} \tau^{1+2(p-1)} \sup_{j \geq 0} \lambda^{t/q} (P^j |f|^q)$$

and

$$c_2 = \left[ \sum_{j=0}^{\infty} j^{-wr} \right]^{t/r} \left[ \sum_{j=0}^{\infty} j^{-\eta-1} \right]^{t/q}. \qquad \square$$

## APPENDIX B

PROOF OF LEMMA 3.2. Let $H_j = P_j - s_j \otimes \nu_j$ for $j = 1, 2$. We begin by showing that

(B.1) $$\{W_{\tau_k^1}\} \overset{\mathrm{d}}{=} \{\underset{\sim}{W}_k\} \qquad \text{when } \underset{\sim}{\lambda} = \widetilde{\lambda},$$

where

(B.2) $$\widetilde{\lambda} = \lambda_2 \Phi_{\lambda_1}, \qquad \Phi_{\lambda_1} = \sum_{\ell=0}^{\infty} \{\lambda_1 H_1^\ell s_1\} P_2^\ell, \qquad \underset{\sim}{P} = P_2 \Phi_{\nu_1}.$$

In order to prove (B.1), it is enough to show that for all integers $r \geq 0$ and for all $A_i \in \mathcal{E}_2^+$,

(B.3) $$\mathbb{P}_\lambda(W_{\tau_0^1} \in A_0, \dots, W_{\tau_r^1} \in A_r) = \mathbb{P}_{\widetilde{\lambda}}(\underset{\sim}{W}_0 \in A_0, \dots, \underset{\sim}{W}_r \in A_r).$$

Let $k_0 = j_0$ and $k_\ell = j_0 + j_1 + \cdots + j_\ell$ for $\ell = 0, \dots, r$. We have

$$\mathbb{P}_\lambda(W_{\tau_0^1} \in A_0 \cdots W_{\tau_r^1} \in A_r)$$

$$= \sum_{j_0=0}^{\infty} \sum_{j_1=1}^{\infty} \cdots \sum_{j_r=1}^{\infty} \mathbb{P}_{\lambda_2}(W_{k_0} \in A_0 \cdots W_{k_r} \in A_r) \mathbb{P}_{\lambda_1}(\tau_0^1 = j_0, \dots, \tau_r^1 = j_r)$$

$$= \sum_{j_0=0}^{\infty} \sum_{j_1=1}^{\infty} \cdots \sum_{j_r=1}^{\infty} \{\lambda_2 P_2^{j_0} I_{A_0} \cdots P_2^{j_r} I_{A_r} 1\} \{\{\lambda_1 H_1^{j_0} s_1\} b_{j_1} \cdots b_{j_r}\}$$

$$= \widetilde{\lambda} I_{A_0} \underset{\sim}{P} I_{A_1} \cdots \underset{\sim}{P} I_{A_r} 1,$$

where $b_\ell = \nu_1 H_1^{\ell-1} s_1$, $\ell \geq 0$. Hence, (B.3) holds.

From

$$\underset{\sim}{P} \geq (s_2 \otimes \nu_2)\Phi_{\nu_1} = s_2 \otimes \nu_2 \Phi_{\nu_1},$$

we obtain the minorization inequality (3.12). Let $\underset{\sim}{H} = \underset{\sim}{P} - \underset{\sim}{s} \otimes \underset{\sim}{\nu}$. Then

$$\underset{\sim}{H} = P_2 \Phi_{\nu_1} - s_2 \otimes \nu_2 \Phi_{\nu_1} = (P_2 - s_2 \otimes \nu_2)\Phi_{\nu_1} = H_2 \Phi_{\nu_1}, \qquad \underset{\sim}{Q} = Q_2 \Phi_{\nu_1},$$

where $Q_2$ (in terms of $H_2$, $s_2$) and $\underset{\sim}{Q}$ (in terms of $\underset{\sim}{H}$, $\underset{\sim}{s}$) are defined by (2.3). The next task is to prove (3.13), that is,

$$\{\widehat{W}_{\tau_k^1}\} \overset{\mathrm{d}}{=} \{\widehat{W}_k\} \qquad \text{when } \underset{\sim}{\lambda} = \widetilde{\lambda},$$

where $\{\widehat{W}_k\}$ denotes the split chain generated by $\underset{\sim}{P}$ and $(\underset{\sim}{s}, \underset{\sim}{\nu})$. Let $\widehat{\underset{\sim}{P}}$ denote the transition probability function for this split chain and let $\widehat{\underset{\sim}{P}}'$ denote the transition probability for $\{\widehat{W}_{\tau_k^1}\}$. We must prove that

$$\widehat{\underset{\sim}{P}}' = \widehat{\underset{\sim}{P}}.$$

First, we recall the structure of a split chain. Suppose that $P$ is a transition probability which satisfies $P \geq s \otimes \nu$. Then the corresponding split chain has transition probability $\widehat{P}$, which satisfies, for $n \geq 1$,

$$\widehat{P}^n(x_0 \times y_0, dx \times y) = y_0 \nu P^{n-1}(dx)\{ys(x) + (1-y)(1-s(x))\}$$
$$+ (1-y_0)QP^{n-1}(x_0, dx)\{ys(x) + (1-y)(1-s(x))\}.$$

In our case, this gives, for $n = 1$,

$$
\begin{aligned}
\text{(B.4)} \qquad \widehat{\underset{\sim}{P}}(w_0 \times y_0, dw \times y) &= y_0 \underset{\sim}{\nu}(dw)\{y\underset{\sim}{s}(w) + (1-y)(1-\underset{\sim}{s}(w))\} \\
&\quad + (1-y_0)\underset{\sim}{Q}(w_0, dw) \\
&\quad \times \{y\underset{\sim}{s}(w) + (1-y)(1-\underset{\sim}{s}(w))\}.
\end{aligned}
$$

We more carefully consider $\widehat{\underset{\sim}{P}}'$, which by (B.2) satisfies

$$\widehat{\underset{\sim}{P}}' = \sum_{\ell=1}^{\infty} b_\ell \widehat{P}_2^\ell.$$

We replace $\widehat{P}_2^\ell$ by the right-hand side of the expression

$$
\begin{aligned}
\widehat{P}_2^\ell(w_0 \times y_0, dw \times y) &= y_0 \nu_2 P_2^{\ell-1}(dw)\{ys_2(w) + (1-y)(1-s_2(w))\} \\
&\quad + (1-y_0)Q_2 P_2^{\ell-1}(w_0, dw) \\
&\quad \times \{ys_2(w) + (1-y)(1-s_2(w))\},
\end{aligned}
$$

where

$$\text{(B.5)} \qquad \sum_{\ell=1}^{\infty} b_\ell \nu_2 P_2^{\ell-1} = \nu_2 \Phi_{\nu_1} = \underset{\sim}{\nu}, \qquad \sum_{\ell=1}^{\infty} b_\ell Q_2 P_2^{\ell-1} = Q_2 \Phi_{\nu_1} = \underset{\sim}{Q}.$$

We then obtain (3.13) from (B.4)–(B.5). The first equality in (3.14) follows from (3.13) and the second is the occupation formula given by (2.10). Actually, when $\lambda = \lambda_1 \times \pi_2$, we get

$$\widetilde{\lambda} = \pi_2 \Phi_{\lambda_1} = \sum_{\ell=0}^{\infty} \{\lambda_1 H_1^\ell s_1\} \pi_2 P^\ell = \pi_2 \{\lambda_1 G_{s_1,\nu_1} s_1\} = \pi_2.$$

Finally, if $\lambda = \nu = \nu_1 \times \nu_2$, then $\widetilde{\lambda} = \underset{\sim}{\nu}$ and $\underset{\sim}{\nu} G_{\underset{\sim}{s},\underset{\sim}{\nu}} = \pi_2$ since $\underset{\sim}{\pi}_{\underset{\sim}{s}} = \pi_{s_2}$.  $\square$

PROOF OF LEMMA 3.3. The waiting times $\{\delta_j, \, j \geq 0\}$ are given by $\delta_j = \tau_j^1 - \tau_{j-1}^1$. Let $b_{n,k} = P(\delta_1 + \cdots + \delta_n = k)$ for $k \geq n$ and $b_{1,k} = b_k$. Then

$$b_{n,k} = \begin{cases} \nu_1 H_1^{k-1} s_1, & \text{for } n = 1 \text{ and } k \geq 1, \\ b_k^{\star n}, & \text{for } n \geq 1 \text{ and } k \geq 1, \end{cases}$$

where "$\star n$" denotes $n$-times convolution. The $n$-step transition probability $\underset{\sim}{P}_n$ is given by

$$\text{(B.6)} \qquad \underset{\sim}{P}_n = \sum_{j=0}^{\infty} b_{n,n+j} P_2^{n+j}, \qquad n \geq 1.$$

Since $\{W_t\}$ is geometric ergodic ([24], Theorem 6.14, page 120), there exist a nonnegative function $M$ such that $\pi_2(M) < \infty$ and a constant $\rho \in (0,1)$ such that

$$\|P_2^n(x,\cdot) - \pi_2\| \leq M(x)\rho^n, \qquad x \in E, \, n \geq 0.$$

Thus, by (B.6),

$$
\begin{aligned}
\text{(B.7)} \qquad \|\underset{\sim}{P}_n(x,\cdot) - \pi\| &\leq \sum_{j=0}^{\infty} b_{n,n+j} \|P_2^{n+j}(x,\cdot) - \pi\| \\
&\leq \sum_{j=0}^{\infty} b_{n,n+j} M(x)\rho^{n+j} \\
&\leq M(x)\rho^n \sum_{j=0}^{\infty} b_{n,n+j}\rho^j \\
&\leq M(x)\rho^n.
\end{aligned}
$$

Hence, by (B.7), $\{\underset{\sim}{W}_k\}$ is geometric ergodic.

For the ergodic $\{W_t\}$, we have

$$\alpha_\ell = \sup_{A,B \in E} \theta_\ell(A, B), \qquad \theta_\ell(A, B) = \pi_2 I_A P_2^\ell I_B 1 - \pi_2 1_A \pi_2 1_B.$$

Here,

$$\theta_\ell(A, B) = \pi_2 I_A P_\ell I_B 1 - \pi_2 1_A \pi_2 1_B = \sum_{j=\ell}^\infty b_{\ell,j} \{ I_A P_2^j I_B - \pi_2 1_A \pi_2 1_B \}$$

$$= \sum_{j=\ell}^\infty b_{\ell,j} \theta_j(A, B).$$

That is,

$$(B.8) \qquad \underset{\sim}{\alpha}_\ell \le \sum_{j=\ell}^\infty b_{\ell,j} \sup_{A,B \in E} \theta_j(A, B) = \sum_{j=\ell}^\infty b_{\ell,j} \alpha_j \le \alpha_\ell.$$

By [5], in general,

$$\sum_{\ell=1}^\infty \ell^k \alpha_\ell < \infty \quad \Longrightarrow \quad \mathbb{E}_{\pi_2} \tau_0^{k+1} < \infty.$$

By (B.8) it follows that

$$\sum_{\ell=1}^\infty \ell^k \alpha_\ell < \infty \quad \Longrightarrow \quad \sum_{\ell=1}^\infty \ell^k \underset{\sim}{\alpha}_\ell < \infty.$$

Hence, (3.15) is true. $\square$

REMARK B.1. We see that $\underset{\sim}{\alpha}_\ell \le \mathbb{E}[\alpha(\delta_1 + \cdots + \delta_\ell)]$. A sharper inequality would be $\underset{\sim}{\alpha}_\ell \le \alpha_{\ell^{1/\beta}} + \mathcal{O}(1)$, and if this inequality is correct, then

$$\sum_{\ell=1}^\infty \ell^{\beta k} \alpha_\ell < \infty \quad \Longrightarrow \quad \sum_{\ell=1}^\infty \ell^k \underset{\sim}{\alpha}_\ell < \infty.$$

**Acknowledgments.** Part of this paper was written while the first author was visiting the Department of Statistics, Stanford University from September of 2003 to August 2004. We are grateful to an Associate Editor and to the referees for a number of valuable comments, not the least for the queries leading to the writing of Section 4.

H. A. KARLSEN
D. TJØSTHEIM
DEPARTMENT OF MATHEMATICS
UNIVERSITY OF BERGEN
JOHS. BRUNSGATE 12
5008 BERGEN
NORWAY
E-MAIL: karlsen@mi.uib.no
         dag.tjostheim@mi.uib.no

T. MYKLEBUST
UNIVERSITY COLLEGE OF SOGN OG FJORDANE
PO BOX 133
6851 SOGNDAL
NORWAY
E-MAIL: terje.myklebust@hisf.no